

\input amstex
\documentstyle{amsppt}
\define\z{{\bold Z }}
\define\cc{{\bold C}}
\define\r{{\bold R}}

\define\q{{\bold Q}}

\define\Ker{\text{Ker}\ }

\define\Pic{\text{Pic}\ }
\define\Tor{\text{Tor}\ }
\define\tor{\text{Tor}}
\define\id{\text{id}}
\define\rk{\text{rk}\ }

\magnification1200

\topmatter

\title
ON BRAUER GROUPS OF REAL ENRIQUES SURFACES
\endtitle

\author
V. V. Nikulin and R. Sujatha
\endauthor

\address
Steklov Mathematical Institute,
ul.Vavilova 42, Moscow 117966, GSP-1, Russia
\endaddress

\email
slava\@nikulin.mian.su
\endemail

\address
School of Mathematics, TATA Institute of Fundamental Research,
Homi Bhabha Road, Bombay-400-005, India
\endaddress

\email
sujatha\@tifrvax.bitnet
\endemail

\abstract
Let $Y$ be a real Enriques surface, $_2Br(Y)$ the subgroup of elements of
order 2 of $Br(Y)$, and  $s,s_{or}$, and  $s_{nor}$
the number of all connected,
connected orientable, and connected non-orientable components of
$Y(\r)$ respectively. Using universal covering $K3$-surface
$X$ of $Y$, we connect $\dim _2Br(Y)$ with
$s, s_{or}$ and $s_{nor}$. As a geometric corollary of our
considerations, we show that $s\le 6$ and $s_{nor}\le 4$.
\endabstract

\endtopmatter

\document

\head
\S 0. Introduction
\endhead

Let $X$ be a smooth projective irreducible  real surface
and
$$
Br(X)=H^2_{et}(X; {\bold G}_m)
$$
denote the Brauer group (cohomological)
of $X$.
This group is known to be a
birational invariant of the surface. Let $X(\r )$ denote the space of
$\bold R$-rational points of $X$ with the Euclidean topology and $s$ be
the number
of real connected components of this space.
Let $_2Br(X)$ be the group of 2-torsion elements in $Br(X).$ If $P\in X$ is
a real point of $X$, we get a natural map
$_2Br(X)\to \  _2Br(P)\cong \z /2$.
In \cite{CT-P}, it is proved that this map
$_2Br(X)\to  \z /2$ does not depend on a choice of the point $P$ in
a connected component of $X(\r )$. Thus, the canonical
 map
$$
_2Br(X)\to (\z /2)^s
\tag0--1
$$
is defined where $s$ is the number
of connected components of $X(\r )$.

A calculation of $\dim \ _2Br(X)$ and the dimension of the kernel of the
map (0--1) is a very interesting and difficult problem. It is important for
a calculation of a such invariant of $X$ as the Witt group $W(X)$, see
\cite{Su}. Probably, these calculations are now known only for real
surfaces  $X$ rational over $\cc $. In this case \cite{Su}, if $s>0$, the
dimension $\dim \ _2Br(X)=2s-1$, the map (0--1) is epimorphic
\cite{CT-P, 3.1}, and
$W(X)\cong (\z)^s\oplus (\z/2)^{s-1}$.
We mention that in \cite{CT-P, 3.1}
it is proved that the map (0--1) is an epimorphism for surfaces $X$
with $H^3(X(\cc );\z/2)=0$. In particular, it is true for rational surfaces.

In this paper, we start studying similar problems for real Enriques surfaces.
Here, for a real Enriques surface $Y$, we calculate
$\dim \ _2Br(Y)$ and connect it with the numbers $s$, $s_{or}$ and
$s_{nor}$ of
connected, orientable and non-orientable components of
$Y(\r )$ respectively.

By a real Enriques surface $Y$,
we mean that $Y_\cc =Y\otimes \cc$ is an irreducible
non-singular
minimal projective algebraic
surface over $\bold C$ with invariants
$\varkappa (Y_\cc)=p_g(Y_\cc )=q(Y_\cc )=0$, \cite{A}.
For Enriques surfaces,
the invariants $p_g(Y_\cc )$ and $q(Y_\cc )$ are the same as for
rational surfaces. But $H^3(Y(\cc );\z/2)=\z/2$.
Thus, it is very interesting to calculate $\dim ~_2Br(Y)$,
study the map (0--1) and $W(Y)$ for real Enriques surfaces $Y$,
and compare these calculations with calculations
for rational surfaces we mentioned above. We mention that
topologically Enriques
surfaces are much more complicated than rational ones
because they are not simply-connected: the fundamental group
$\pi_1(Y(\cc ))\cong \z/2$. It makes these
calculations very delicate.

Now, we formulate basic results of the paper.
Let $Y$ be a real Enriques surface.
Let $Y(\cc )$ denote the underlying complex manifold of $Y_\cc $ and
$G$ be the Galois group $G(\cc  /\r )$. We identify the generator of
this group with the corresponding antiholomorphic involution $\theta $ on
$Y(\cc )$.
Everywhere below for a 2-elementary group
$\dim  (\z /2)^a=a$.
 For any $G$-module $A$, $A^G=A^\theta $
denotes the set of elements of $A$ fixed by $G$.

We introduce the following basic invariants of a real
Enriques surface $Y$:

\definition{Definition 0.1}
The invariant
$\epsilon(Y)=1$ if the differential
$$
d_2^{0,2}:E_2^{0,2}=H^0(Y(\cc );\z /2)^G \to
E_2^{2,1}=H^2(G;H^1(Y(\cc );\z /2))\cong \z/2
$$
of the Hochschild--Serre spectral
sequence (see \cite{Mi}) vanishes,
and $\epsilon(Y)=0$ otherwise.
For Enriques surfaces, $H^1(Y(\cc );\z /2)\cong \z/2$.

The invariant
$$
b(Y) = \dim H^2(Y(\cc );\z /2)^\theta -
\dim (\Pic Y_\cc )^\theta /2(\Pic Y_\cc)^\theta + 1.
$$
\enddefinition

First of all we prove

\proclaim{Theorem 0.1}  Let $Y$ be a real Enriques surface and
$Y(\r )\not= \emptyset $.

Then
$$
\dim~_2Br(Y)=b(Y) + \epsilon (Y).
$$
\endproclaim

\demo{Proof} See Theorem 1.2 (in \S  1).
\enddemo

Theorem 0.1 shows that to estimate $\dim~_2Br(Y)$, it is important to
estimate the invariant $b(Y)$,
because $0\le \epsilon (Y)\le 1$.

We prove the following general estimates from below
for the $b(Y)$ and
\linebreak
$\dim~_2Br(Y)$.

\proclaim{Theorem 0.2} Let $Y$ be an arbitrary real Enriques surface
and $s$ the number of connected components of $Y(\r )$. Then
$$
b(Y)\ge 2s-2.
$$

It follows (by Theorem 0.1)
$$
\dim~_2Br(Y)\ge 2s-2 + \epsilon (Y)\ge 2s-2.
$$
\endproclaim

\demo{Proof} See Theorem 2.1 (in \S  2).
\enddemo

Considering the image of the $_2Br(\text{Spec}\  \r)$ of the basic field
$\r $ in
$_2Br(X)$, we evidently get

\proclaim{Proposition 0.3}  For any real surface $X$,
the element
$$
(\undersetbrace \text{$s$ times} \to{1, \cdots ,1})
$$
belongs to the
image of the map (0--1).
\endproclaim

Thus, from these three statements we get the

\proclaim{Corollary 0.4} For a real Enriques surface $Y$
$$
\dim_2Br(Y)\ge s.
$$
\endproclaim

This corollary is interesting because it is not now known that the
map (0--1) is an epimorphism for real Enriques surfaces.

In \S 3 of the paper, we give some precise formula for the invariant
$b(Y)$. This is the most non-trivial result of the paper. We recall
\cite{C-D} that the universal covering surface of an Enriques surface
$Y$ is a $K3$-surface $X$, and the universal covering
$\pi :X(\cc ) \to Y(\cc )$ is 2-sheeted.
We denote by $\tau$ the holomorphic involution of $X$ corresponding to
this covering. It is not difficult to see that if
$Y(\r )\not= \emptyset$, then there are two liftings $\sigma $ and
$\tau\sigma $ of the antiholomorphic involution $\theta $ of $Y(\cc )$ to
antiholomorphic involutions of $X(\cc )$.

\proclaim{Theorem 0.5}
Let $Y$ be a real Enriques surface and suppose there are
two liftings $\sigma $ and $\tau \sigma $
of antiholomorphic involution $\theta $ of
$Y(\cc )$ to antiholomorphic
involutions of the universal covering $K3$-surface $X(\cc )$
(for example this
is true if $Y(\r )\not= \emptyset $). Let $s_{or}$ and $s_{nor}$ be the
number of orientable and non-orientable connected
components of $Y(\r )$ respectively.

Then
$$
b(Y)=s_{nor}+2s_{or} - z(\sigma)+ \dim H(\sigma )_- -
\dim  (H(\sigma )_+)^\perp \cap H(\sigma )_- + \beta (Y).
$$
\endproclaim

\demo{Proof} See Theorems 3.5.1, 3.5.2 and 3.5.3 and the formula
 (3--5--1) (in \S 3).
\enddemo

Here we cannot give precise definitions of the summands of this
formula. We only mention the following inequalities for these integers:
$0\le z(\sigma)\le 2$,
$0 \le \dim H(\sigma )_- - \dim  (H(\sigma )_+)^\perp \cap H(\sigma )_-\le 2$,
$0 \le \beta (Y)\le 1$.
See \S 3 for precise definitions of these
integers. It is very important that integers
$s_{nor}+2s_{or}$,  $z(\sigma)$,
$\dim H(\sigma )_-$ and
$\dim (H(\sigma )_+)^\perp \cap H(\sigma )_-$ are
defined using only the action of the
involutions $\sigma $ and $\tau \sigma $ on $H^2(X(\cc );\z)$. In some
cases we can prove that the invariant $\beta (Y)=0$. Thus, using global
Torelli Theorem \cite{P\u S-\u S}
and epimorphicity of Torelli map \cite{Ku} for $K3$-surfaces, and results of
\cite{N4},  we can construct
Enriques surfaces with  these prescribed  invariants.

By the inequalities above, Theorem 0.5 gives the estimate for $b(Y)$ from
above. Since $2s=2s_{nor}+2s_{or}$ and
Theorem 0.2 gives the inequality for $b(Y)$ from below,
the formula of Theorem 0.5 and the inequality of Theorem 0.2 give
very strong estimates
on the numbers  $s_{nor}, s_{or}$ and show that the inequality of the
Theorem 0.2 is not far from being an equality.  For a classification of real
Enriques surfaces using these considerations,
see \cite{N6}.

We mention that in Theorem 3.4.7 we give another formula  for the number
$b(Y)$ which is more useful in some cases.

 From Theorem 0.5 and statements 0.1 ---0.3 above,
we get the following result when we calculate
$\dim~_2Br(Y)$ precisely. To formulate this result, we should introduce
some invariants of $Y$.

Let $H^2=H^2(Y(\cc );\z)/\text{Tor}$.
The $H^2$ is an unimodular lattice with
respect to the intersection pairing. The involution $\theta$ acts on $H^2$.
The invariant $r(\theta)$ is defined by $r(\theta)=\rk (H^2)^\theta$.
The intersection pairing defines a group
$A_{(H^2)^\theta} = ((H^2)^\theta)^\ast/(H^2)^\theta$. Since the lattice
$H^2$ is unimodular, this group is $2$-elementary:
$A_{(H^2)^\theta}\cong (\z/2)^{a(\theta)}$. This defines the invariant
$a(\theta)$. Another definition of this invariant gives the formula
$\dim~(H^2/2H^2)^\theta = \rk H^2-a(\theta)$.

\proclaim{Theorem 0.6} Let $Y$ be a real Enriques surface and
$Y(\r )\not= \emptyset $.

Then
$b(Y)=0$
iff the surface $Y(\r )$ is connected non-orientable and for the
invariants $r(\theta)$ and $a(\theta )$ we have the equality
$r(\theta )=a(\theta )$.

By Theorem 0.1 and Proposition
0.3, for this surface $Y$ the
invariant $\epsilon (Y)=1$ and $\dim~_2Br(Y)=1$.
\endproclaim

\demo{Proof} See Theorem 3.6.1.
\enddemo

The idea of the proof of Theorem 0.1 is  to use the Kummer sequence
and estimate the dimension of $_2Br (Y)$ using the Hochschild--Serre
spectral sequence. The proof of Theorem 0.2 is not difficult and uses
Lefschetz formula for fixed points and Smith exact sequence.
The proof of the
Theorems 0.5 and 0.6 is hard and uses results of
\cite{Ha1}, \cite{N3}, \cite{N4}.
  The most important ingredient of our proof is the use of the
theory of involutions of lattices (integral quadratic forms)
with condition on sublattice which was developed in
\cite{N4} (and also \cite{N3}).
We apply this theory to the action of antiholomorphic involutions
on the 2-cohomology lattice of $Y(\cc )$ and  its universal covering
$K3$-surface $X(\cc )$. The proof of Theorems 0.2, 0.5 and 0.6
is a part of the general problem of
the classification of real Enriques surfaces which is studied in \cite{N6}
(as an example, see Sect. 3.7 at the end of the paper).

All the  results of this paper may be generalized to the following
more general situation: Instead of $K3$-surfaces one should consider
 complex smooth projective algebraic surfaces $X$
such that 2-torsion $_2\Pic X=0$. Instead of Enriques surfaces one should
consider real surfaces $Y=X/\{ id,\tau \} $ where $\tau $ is a holomorphic
involution of $X$ without fixed points. For example, statements 0.1 -- 0.4
are true in this case. The situation is more complicated with
Theorems 0.5 and 0.6 (see
\S 3.4 and Lemma 3.4.3) but it is similar.
 We hope to generalize the results here for this more general case
in subsequent publications.

In Sect. 3.8, we cite further results on real Enriques surfaces
from \cite{N5} and \cite{N6}, which were obtained by the
first author during the time this paper was considered for publication.

We are grateful to Prof. J--L. Colliot-Th\'el\`ene for useful comments.

The first author worked on this article during his stay at the first
half of 1991 in  Tata Institute of Fundamental Research, Bombay;
Steklov Mathematical
Institute, Moscow and Institut des Hautes \'Etudes Scientifiques, Paris.
He is grateful to these Institutes for hospitality. Preliminary
variant of this paper was published as a preprint \cite{N-S}.

\newpage

\head
\S  1. Estimation of the 2-torsion of
Brauer group of\\
real Enriques surfaces
\endhead

At first, we give some basic facts which we will use
about complex Enriques surfaces.

We recall that a complex $K3$-surface is a smooth projective algebraic
surface $X$ over $\cc$ such that $X(\cc )$ is a simply-connected complex
surface and the canonical class $K_X=0$. It follows that there exists a
2-dimensional regular differential form $\omega_X\in \Omega [X]$ such
that the divisor $(\omega_X)=0$. Thus, for an arbitrary
point $z\in X(\cc )$ and local coordinates $z_1,z_2$ at a neighbourhood of
$z$, we can
write down this form as $f(z_1,z_2)dz_1\wedge dz_2$ where
$f(z_1,z_2)$ is a holomorphic function and
$f(z)\not=0$. This form $\omega_X$ is unique up to multiplication on elements
of $\cc $.

We recall basic facts about $K3$-surfaces $X$ we will use
(see G. N. Tjurina,
\cite{A, Ch. 9}).
We have:
$$
\split
& H^0(X(\cc );\z)\cong H^4(X(\cc );\z)\cong\z,\\
& H^1(X(\cc );\z)\cong H^3(X(\cc );\z)\cong 0,\\
& H^2(X(\cc );\z)\cong \z^{22}.
\endsplit
\tag1--0
$$
Besides, the group $H^2(X(\cc );\z)$ with the intersection
pairing is an even unimodular lattice of signature $(3,19)$.

By definition, an Enriques surface
$Y$ over $\cc $ is a minimal smooth projective
algebraic surface over $\cc $ with invariants
$\varkappa (Y)=p_g(Y)=q(Y)=0$ (see B. G. Averbuh, \cite{A, Ch. 10}).
An equivalent definition of Enriques surfaces is that
an Enriques surface $Y$ is
the quotient $Y=X/\{ \id, \tau \}$ of a $K3$-surface
$X$ by an algebraic involution
$\tau$ without fixed points.
In the paper, we will use the last
definition of Enriques surfaces. Thus, $X(\cc )$ is
the universal covering of $Y(\cc )$, and
$$
\pi_1(Y(\cc ))=H_1(Y(\cc );\z)\cong \z/2.
\tag1--1
$$
By \cite{N1, \S 5},
$$
\tau^*(\omega_X)=-\omega_X.
\tag1--2
$$
for an involution $\tau $ without fixed points on a $K3$-surface.
It follows that the canonical class of $Y$
$$
K_Y\not=0,\  but\  2K_Y=0.
\tag1--3
$$
Using (1--0), (1--1) and
standard topological facts:
Lefschetz fixed point formula, universal coefficient formula,
Poincar\'e duality \cite{Sp}
and elementary
facts about actions of finite groups \cite{Br, Ch. II, Sects. 2, 3}, one
can prove easily that:
$$
\split
& H^0(Y(\cc );\z)\cong \z,\  H^1(Y(\cc );\z)=0,\\
& H^2(Y(\cc );\z)\cong \z^{10}\oplus \z/2,\\
& H^3(Y(\cc );\z)\cong \z/2,\  H^4(Y(\cc );\z)\cong \z;
\endsplit
\tag1--4
$$
$$
\split
& H^0(Y(\cc );\z/2)\cong \z/2,\  H^1(Y(\cc );\z/2)\cong \z/2,\\
& H^2(Y(\cc );\z/2)\cong (\z/2)^{12},\\
& H^3(Y(\cc );\z/2)\cong \z/2,\  H^4(Y(\cc );\z/2)\cong \z/2;
\endsplit
\tag1--5
$$
and we have an isomorphism:
$$
\pi^\ast : H^2(Y(\cc );\z)/\tor \to H^2(X(\cc );\z)^\theta
\tag1--6
$$
where $\pi :X \to Y$ is the quotient morphism.
Besides,
since $H^1(Y(\cc );\z)=0$, the characteristic class map gives an isomorphism
$\Pic Y \cong H^2(Y(\cc );\z )$. By Poincar\'e duality and
Hodge index theorem,
the $H^2(Y(\cc );\z)/\tor $
with intersection pairing is an
unimodular lattice of signature $(1,9)$. By the formula for the genus of
a curve on an algebraic surface and (1--2), this lattice is even.

\smallpagebreak

Now let us consider real Enriques surfaces.
By definition, a real Enriques surface is a smooth projective algebraic
surface $Y$ over $\r$ such that $Y_\cc=Y\otimes _\r \cc$ is a complex
Enriques surface.

Let $G=Gal(\cc /\r )=\{\id,\theta \}$. The $\theta $ acts on
$Y(\cc )$ as an antiholomorphic involution and
$Y(\r )=Y(\cc )^G$. Here and as follows,
$U^G=U^\theta$ denote the set of fixed elements for an action of
a group $G$ on a set $U$.

We recall \cite{Mi} that the Kummer exact sequence
$$
0 \to \mu_2 \to {\bold  G}_m \to  {\bold  G}_m \to 0
$$
yields the exact sequence
$$
0 \to \Pic Y/2\Pic Y \to H^2_{et}(Y; \mu _2) \to \ _2Br(Y) \to 0.
\tag1--7
$$
Since $Y(\r )$ is non-empty, $\Pic Y=(\Pic Y_\cc)^G$,
(see \cite{Ma}). (One can deduce this from Hochschild--Serre spectral
sequence and
the canonical isomorphisms \linebreak
$\Pic Y\simeq H^1_{et}(Y;{\bold G}_m),\
\Pic Y_\cc \simeq H^1_{et}(Y_\cc;{\bold G}_m)$, see \cite{Mi}.)
By (1--7), we then get
$$
\dim \ _2Br(Y)=\dim~H^2_{et}(Y;\mu_2)-
\dim (\Pic Y_\cc)^G/2(\Pic Y_\cc)^G.
\tag1--8
$$
The dimension of the \'etale cohomology group $H^2_{et}(Y;\mu_2)$ is
estimated using the Hochschild--Serre
spectral sequence
$$
E_2^{p,q}=H^p(G;H_{et}^q(Y_\cc;\mu_2))
\Longrightarrow H^{p+q}_{et}(Y;\mu_2)=E^{p+q}.
$$
where for a complex algebraic manifold $Y_\cc$ we have
$$
H_{et}^q(Y_\cc ;\mu_2)=H^q(Y(\cc );\z/2),
$$
see \cite{Mi}.

The following invariant of a real Enriques surface is very important.
We recall \linebreak
$H^1(Y(\cc );\z /2)\simeq \z /2$, by (1--5).

\definition{Definition 1.1} For a real Enriques surface $Y$ an
invariant
$\epsilon(Y)=1$ if the differential
$d_2^{0,2}:E_2^{0,2}=H^0(Y(\cc );\z /2)^G \to
E_2^{2,1}=H^2(G;H^1(Y(\cc );\z /2))\simeq \z/2$ of the
 Hochschild--Serre spectral
sequence vanishes,
and $\epsilon(Y)=0$ otherwise.
\enddefinition

We then have the following Lemma,

\proclaim{Lemma 1.1} Let $Y$ be a real Enriques surface and
$Y(\r )\not=\emptyset$.
Then
$$
\dim H^2_{et}(Y;\mu_2   )=1+\dim H^2(Y(\cc );\z /2)^G + \epsilon(Y).
$$
\endproclaim

\demo{Proof} Denoting $H^2_{et}(Y; \mu_2 )$ by $E^2$
and using the usual spectral
sequence notation, we see that
$$
\dim~ E^2 = \dim~E_{\infty}^{2,0} + \dim~E_{\infty}^{1,1} +
\dim~E_{\infty}^{0,2}.
\tag 1--9
$$
We calculate
$\dim~E_{\infty}^{2,0}$, $\dim~E_{\infty}^{1,1}$ and
$\dim~E_{\infty}^{0,2}$ below.

Since $Y(\r )$ is not empty, consider a real point
$P \in Y(\r )$.  By functoriality,
this induces a spectral
sequence homomorphism,
$$
H^p(G; H^q(Y(\cc );\z /2)) \to H^p(G; H^q(P;\z /2)).
$$
But the homomorphism
$$
H^p(G;H^0(Y(\cc );\z /2)) \to H^p(G; H^0(P;\z /2))
$$
is an isomorphism and the spectral sequence for the point is trivial.
Therefore in the spectral sequence for $Y$, the differentials $d_r^{p-r,r-1}$
are identically zero. This implies that $E_{\infty}^{2,0} \simeq E_2^{2,0}
\simeq \z /2$, and by (1--5), $E_{\infty}^{1,1} \simeq E_2^{1,1}
\simeq \z/2$. Thus,
$\dim~E_{\infty}^{2,0}=\dim~E_{\infty}^{1,1}=1$.
Further, by the above remark, the differential
$d_3^{0,2}:E_3^{0,2} \to E_3^{3,0} $
is zero, and we have
$E_3^{0,2} \simeq E_{\infty}^{0,2}$
and an exact sequence
$$
0 \to E_{\infty}^{0,2} \to E_2^{0,2} \buildrel {d_2^{0,2}} \over
\longrightarrow E_2^{2,1}
$$
where by (1--5),
$E_2^{2,1}=H^2(G; \z /2) \simeq \z /2$.
Thus, by the definition 1.1,
$$
\dim~E_2^{0,2}=
\dim~E_{\infty}^{0,2} + 1- \epsilon (Y).
\tag 1--10
$$
By definition, $E_2^{0,2} =H^0(G; H^2(Y(\cc ); \z /2)) =
 H^2(Y(\cc );\z /2)^G$, and by (1--10),
 $$
 \dim~E_{\infty}^{0,2}=\dim~H^2(Y(\cc );\z /2)^G - 1+\epsilon (Y).
 $$

This proves the Lemma.
\enddemo

\definition{Definition 1.2} We denote by $b(Y)$ the following
invariant of a real Enriques surface $Y$ with an antiholomorphic
involution $\theta $:
$$
b(Y)=\dim H^2(Y(\cc );\z /2)^\theta
-\dim (\Pic Y_\cc)^\theta/2(\Pic Y_\cc)^\theta +1.
$$
\enddefinition

By Lemma 1.1 and (1--8), we get the following

\proclaim{Theorem 1.2}
Let $Y$ be a real Enriques surface with an antiholomorphic involution
$\theta $ and $Y(\r )\not= \emptyset $.

Then
$$
\dim~_2Br(Y)=b(Y)+\epsilon (Y)
$$
where
$\epsilon (Y)=0$ or $1$. See Definitions 1.1 and 1.2 above.
\endproclaim

Theorem 1.2 shows that to estimate $\dim~_2Br(Y)$ it is very important
to estimate the invariant $b(Y)$.
The following \S 2  and \S 3 will be devoted to this problem. We want to
connect $b(Y)$ with the numbers $s$, $s_{nor}$ and $s_{or}$ of
connected, non-orientable and orientable components of $Y(\r )$.

\newpage

\head
\S 2. Inequalities:  $b(Y)\ge 2s-2,\ \dim~_2Br(Y)\ge 2s-2+\epsilon (Y)$
\endhead

For a real rational surface $X$ with $X(\r )\not= \emptyset$
we have \cite{Su}: $\dim~_2Br(X)=2s-1$.
In this paper, for a real Enriques surface $Y$, we can only prove similar
inequalities.

\proclaim{Theorem 2.1} Let $Y$ be an arbitrary real Enriques surface.

Then $b(Y) \ge 2s-2$.

It follows (by Theorem 1.2) that
$\dim~_2Br(Y)\ge 2s-2+\epsilon (Y)$.
\endproclaim

By Theorem 2.1 and Proposition 0.3, we also get the following result which
is interesting because for a real Enriques surface we don't know that
the homomorphism (0--1) is an epimorphism.

\proclaim{Corollary 2.2} Let $Y$ be an arbitrary real Enriques surface
and $s$ the number of connected components of $Y(\r )$.

Then $\dim~_2Br(Y)\ge s$.
\endproclaim

\demo{Proof of Theorem 2.1} We introduce an important invariant of
a real Enriques surface $Y$ with an
antiholomorphic involution $\theta$:
$$
r(\theta )=\rk (H^2(Y(\cc );\z )/\tor)^\theta .
\tag2--1
$$
We prove the following
formula  for an arbitrary real Enriques surface:
$$
\dim (\Pic Y_\cc )^\theta /2(\Pic Y_\cc )^\theta =
11-r(\theta ).
\tag2--2
$$
By (1--4), $H^1(Y(\cc);\r)=0$. Thus, we have the canonical
isomorphism of characteristic class
$$
\Pic Y_\cc\cong H^2(Y(\cc );\z ).
$$
For this isomorphism, the action of the
involution $\theta $ on the Picard group goes to the action of
the involution $-\theta $
on the $H^2(Y(\cc ); \z )$
(this is well-known, for example, see
\cite{Si, Sect. I, 4}).
It follows therefore, that we have an isomorphism,
$$
(\Pic Y_\cc )^\theta \cong H^2(Y(\cc );\z )_\theta .
$$
Here for a module $M$ and an involution $\theta $ of this module we
denote $M_\theta =\{ x \in M \mid \theta (x)=-x \} $.

By (1--4),  $H^2(Y(\cc );\z )\cong
\z /2\oplus \z ^{10}$ where
$$
\z ^{10}\cong
H^2(Y(\cc );\z )/\Tor H^2(Y(\cc );\z ).
$$
We use the following obvious statement: Let $M=\z /2\oplus \z ^m$
and a group $G$ of order two acts on $M$. Then
$M^G\cong \z /2\oplus \z ^n$ if
$(M/\Tor M)^G\cong \z ^n$.

By definition of the invariant $r(\theta )$,
it follows that
$$
\rk  (H^2(Y(\cc );\z )/Tor)_\theta =10-r(\theta ).
$$
By the remark above, we get the statement (2--2).

By (1--4), Lefschetz's fixed point formula
(see \cite{Sp, Ch. 5, Sect. 7, Theorem 6})
gives for the involution $\theta $
$$
\chi (Y(\r ))=2r(\theta )-8.
$$
It follows that
$$
r(\theta )=(1/2)\chi (Y(\r ))+4.
\tag2--3
$$

 From \cite{Kr, Theorem 2.3}
or from Smith exact sequence (see \cite{Br, Ch.3}), we have
$$
\dim~H^\ast (Y(\r );\z /2)\le
\dim~H^1(G; H^\ast (Y(\cc ); \z /2)).
$$
Here $H^1(G; H^\ast (Y(\cc ); \z /2))=
H^\ast (Y(\cc ); \z /2)^\theta /(1+\theta)H^\ast (Y(\cc ); \z /2)$.
 From the exact sequence
$$
0 \to H^\ast (Y(\cc ); \z /2)^\theta \to H^\ast (Y(\cc ); \z /2) \to
(1+\theta )H^\ast (Y(\cc ); \z /2) \to 0
$$
we get
$\dim~H^1(G;H^\ast (Y(\cc ); \z /2))=
2\dim~H^\ast (Y(\cc ); \z /2)^\theta -\dim H^\ast (Y(\cc ); \z /2)$.
Thus, from the inequality above, we have
$$
\dim~H^\ast (Y(\r );\z /2)\le
2\dim~H^\ast (Y(\cc ); \z /2))^\theta -
\dim~H^\ast (Y(\cc ); \z /2).
$$
By (1--5), it follows that
$$
\dim H^2(Y(\cc );\z /2)^\theta \ge
4+(1/2)\dim H^\ast (Y(\r );\z /2).
\tag2--4
$$
 From (2--2),
$$
\dim H^2(Y(\cc );\z /2)^\theta -
\dim (\Pic Y_\cc)^\theta /2(\Pic Y_\cc)^\theta +1=$$
$$=\dim H^2(Y(\cc );\z /2)^\theta +r(\theta )-10.
\tag2--5
$$
Thus, from (2--3), (2--4) and (2--5), we get that
$$
\split
b(Y)=\dim H^2(Y(\cc );\z /2)^\theta &-
\dim (\Pic Y_\cc)^\theta /2(\Pic Y_\cc)^\theta +1\\
&\ge
(1/2)\dim H^\ast (Y(\r );\z /2)+
(1/2)\chi (Y(\r ))-2.
\endsplit
\tag2--6$$
For a closed connected surface $F$, we have
$\dim H^\ast (F;\z /2)+\chi (F)=4$.
Thus, the right side
of the inequality (2--6) is equal to $2s-2$.

This proves  the Theorem 2.1.
\enddemo

\newpage

\head
\S 3. Calculation of the invariant $b(Y)$
using universal covering $K3$-surface $X$
\endhead

\subhead
3.1. Notation
\endsubhead

We use the notations set up
in the previous \S  1 and \S 2.
Thus, $Y$ is a real Enriques surface, $Y(\r )$ is
the set of real points of $Y$, and $Y(\cc )$ the set of complex points
with the corresponding antiholomorphic involution $\theta $.
We want to obtain in this section a precise formula for the invariant
 $b(Y)$ (see Definition 1.2)
using the universal covering $K3$-surface $X$.

\subhead
3.2. Real part of Enriques surface $Y$
and the universal covering $K3$-surface
\endsubhead
By definition (see the beginning of \S 1), an Enriques surface is a
quotient surface
$Y_\cc=X/ \{\id,\tau \}$ where $X$ over $\cc $ is a complex
$K3$-surface and $\tau $ an algebraic
involution of $X$
without fixed points. Let $\pi :X \to Y_\cc $ be the
quotient morphism. Since $X(\cc )$ is simply-connected,
$\pi :X(\cc )\to Y(\cc )$ is the 2-sheeted universal
covering with the holomorphic involution $\tau $ of the
covering. Thus, $\pi \tau =\pi $.

Let $\theta $ be the antiholomorphic involution of the complex
surface $Y(\cc )$ corresponding to the real surface $Y$.
Let us suppose that
$$
Y(\r )\not= \emptyset
$$
and
${\bar y}\in Y(\r ).$ Let $y\in X(\cc )$ such that
$\pi (y)={\bar y}.$
Since $\pi $ is the universal covering, there exists an antiholomorphic
automorphism $\sigma $ of $X(\cc )$ such that
$\pi \sigma =\theta \pi $ and $\sigma (y)=y.$ It follows therefore that
$\sigma ^2$ is holomorphic, $\sigma ^2(y)=y$, and $\pi \sigma ^2 =\pi$.
Thus, $\sigma ^2$ is an element of the group of the unramified covering
$\pi $. Since $\sigma ^2$ has the fixed point $y$,
$\sigma ^2=id$.
Hence it follows that $\sigma $ is an antiholomorphic involution
of $X(\cc )$ with the fixed point
$y\in X_\sigma(\r )=X(\cc )^{\sigma }$.
 Evidently, the point
$\tau (y)\in X_\sigma (\r )$ too. Here we denote by $X_\sigma $ the real
$K3$-surface defined by the antiholomorphic involution $\sigma $. Thus,
$X_\sigma (\r )$ is the set of real points of $X_\sigma $.

The composition
${\tilde \sigma }=\tau \sigma $ is also an antiholomorphic automorphism of
$X(\cc )$, $\pi \tilde {\sigma }=\theta \pi $, and
$\tilde \sigma (y)=\tau (y)$. It follows that $\tilde {\sigma }^2=id.$
Thus, $\tilde \sigma =\tau \sigma $
is another lifting on $X(\cc )$ of the antiholomorphic involution
$\theta $  of $Y(\cc )$. Thus, $\tau\sigma $ defines another real
$K3$-surface $X_{\tau\sigma}$ with the set of real points
$X_{\tau\sigma}(\r )
=X(\cc )^{\tau \sigma}$.
As in the case of $X_\sigma (\r )$, the real part
$X_{\tau\sigma}(\r )$ is invariant by the action of $\tau $.
It is obvious that there do not exist other liftings of $\theta $
to $X$ except $\sigma $ and $\tau \sigma $.

Other speaking, we lift the group
$G=\{\id ,\theta\}$ of order 2 on $Y$ to the
group
$\Gamma =\{\id ,\tau ,\sigma ,\tau\sigma\}$
of order $4$ on $X$. The group $\Gamma \cong (\z/2)^2$.
It is difficult to connect the action of $G$ on cohomology of $Y$ with
the real part $Y(\r )$. Therefore, we want to connect the action of
$\Gamma $ on cohomology of $X$ with the
$X_\sigma (\r )$, $X_{\tau\sigma }(\r )$ and $Y(\r )$.

The real parts
$X_\sigma (\r )$ and
$X_{\tau\sigma}(\r )$ have no common point, since such a
point would be the fixed point of the involution $\tau $, but $\tau $
has no fixed point in $X(\cc )$.

Let $\omega _X$ be a non-zero 2-dimensional holomorphic
differential form on X (see the beginning of \S 1).
This form is unique up to multiplication by
non-zero elements of
$\cc $ and it has no zeros on $X(\cc )$.
(Roughly speaking, it is a complex volume form.)
By (1--2),
$\tau ^\ast (\omega _X)=-\omega _X.$ We can choose the form $\omega _X$ by
the condition that $\sigma ^\ast (\omega _X)=\overline {\omega _X}.$ This
defines the form $\omega _X$ up to multiplication by real numbers.
We denote this form as $\omega _X^\sigma $. The form
$\omega _X^\sigma $ gives the canonical
(up to multiplication by non-zero reals) volume form on
the real part $X_\sigma(\r )$ and defines the canonical
orientation on it (up to a change of this orientation on all
connected components of $X_\sigma(\r )$ simultaneously).
See \cite{N3, Theorem 3.10.6} for details. Since
$\tau ^\ast (\omega _X^\sigma )=-\omega _X^\sigma$, the involution
$\tau $ restricted on $X_\sigma(\r ) $
changes this orientation. The same statement holds for the real part
$X_{\tau\sigma}(\r )$. Thus, we have proved the following
very important

\proclaim{Lemma 3.2.1}
Let $Y$ be a real Enriques surface,
$\theta $ the corresponding antiholomorphic
involution on $Y(\cc )$ and the
real part $Y(\r )\not= \emptyset $. Let $\pi :X(\cc )\to
Y(\cc )$ be the universal two sheeted covering where
$X(\cc )$ is $K3$-surface, and $\tau $ the involution
of the covering ${\pi }$.

There are precisely two liftings $\sigma $ and
$\tilde \sigma =\tau \sigma $
of $\theta $ to  antiholomorphic involutions of
$X(\cc )$. Let
$X_\sigma(\r ) =X(\cc )^\sigma $ and
$X_{\tau\sigma}(\r )=X(\cc )^{\tau \sigma }$
be the real parts
corresponding to these antiholomorphic involutions equipped with
canonical orientations defined by the canonical volume forms
$\omega _X^\sigma $ and $\omega _X^{\tau \sigma }$
 respectively. Then
$X_\sigma(\r ) \cap X_{\tau\sigma}(\r )=\emptyset $,
$\tau (X_\sigma(\r ) )=X_\sigma(\r )$ and
$\tau (X_{\tau\sigma}(\r ))=X_{\tau\sigma}(\r )$, and
the involution $\tau $ changes
the canonical orientation of this sets (defined by volume forms
$\omega _X^\sigma $ and $\omega _X^{\tau \sigma }$ respectively)
to the opposite
one. Moreover, with respect to this action of $\tau $ the morphism
$\pi $ is the factorisation morphism which defines the canonical
identification

$$
Y(\r )=X_\sigma(\r ) /\{ id, \tau \} \coprod
X_{\tau\sigma}(\r )/\{ id, \tau \} .
$$
In particular, for the numbers  $s(\sigma )$ and
$s(\tau \sigma )$  of
connected components of $X_\sigma(\r ) $ and
$X_{\tau\sigma}(\r )$ respectively, and the  numbers
$s_{or}$ and $s_{nor}$ of orientable and non-orientable
connected components of $Y(\r )$ respectively, we have the
equality
$$
s(\sigma )+s(\tau \sigma )=s_{nor}+2s_{or}.
$$
\endproclaim

In the \S 2, we connected the invariant $b(Y)$ with the number
$s=s_{or}+s_{nor}$.
In the following sections, we want to connect the $b(Y)$ with the
number $s_{nor}+2s_{or}$.
Lemma 3.2.1 is important for us, because, by this Lemma,
the number $s_{nor}+2s_{or}$
is the same as $s(\sigma )+s(\tau \sigma )$.

\subhead
3.3. Invariants of the action of
$\Gamma =\{\id ,\tau, \sigma , \tau \sigma \}$ on $H^2(X(\cc );\z)$ and
the sets $X_\sigma(\r)$ and $X_{\tau\sigma}(\r)$
\endsubhead

Henceforth, let
$Y$ be a real Enriques surface such that the antiholomorphic
involution $\theta $ of $Y(\cc )$
 has two liftings $\sigma $ and
$\tau \sigma $ to antiholomorphic involutions of
the $K3$-surface $X(\cc )$. By  Lemma 3.2.1,
this is true if $Y(\r )\not= \emptyset $.

Let $L$ be a lattice $H^2(X(\cc );\z )$ with the intersection pairing.
The $L$ is an even unimodular lattice of
signature (3,19) (see the beginning of \S 1).

For the module $M$ and an involution $\phi $ of this module we denote
$$
M^\phi =\{ x\in M \mid \phi (x)=x \},\ \ \ \
M_\phi =\{ x\in M \mid \phi (x)=-x \}.
$$
Let $\phi $ be an involution on a complex $K3$-surface. We suppose that
$\phi$ is either holomorphic
with the condition $\phi^\ast \omega_X=-\omega_X$ or antiholomorphic.
This involution $\phi$ acts on $L$. As for
any involution of an even unimodular lattice,
invariants $$(r(\phi ),a(\phi ),\delta (\phi ))$$ of this action
are defined (see \cite{N3, \S 3}).
Here
$$
r(\phi )=\rk L^\phi .
\tag3--3--1
$$
It is known (see \cite{N2, \S 4} for the holomorphic
case and
\cite{H, appendix}, \cite{N3, \S  3}
for the antiholomorphic one)
that the lattice $L^\phi $ has the
signature $(1,r(\phi )-1)$, i.e. it is hyperbolic.

Since the lattice $L$ is unimodular,
the lattice $L^\phi $ is 2-elementary. This means that
the discriminant group $A_{L^\phi }=(L^\phi)^\ast /L^\phi $
is isomorphic to a 2 -elementary group:
$$
A_{L^\phi }=(L^\phi)^\ast /L^\phi \simeq
 (\z /2)^{a(\phi )}.
\tag3--3--2
$$

(3--3--3) The invariant $\delta (\phi )$ is
equal to 0 if the discriminant quadratic form $q_{L^\phi }$ of
the lattice $L^\phi $ is even, i.e. $q_{L^\phi }(u)\in \z /2\z$
for any $u\in A_{L^\phi }$, and $\delta (\phi )=1$ otherwise.

Here and as follows we denote by $A_S$ the discriminant group
$A_S=S^\ast/S$ of a lattice $S$ and by $q_S$ and $b_S$ the
discriminant quadratic and bilinear form of the lattice $S$ (the
discriminant quadratic form is defined for an even lattice $S$ only).
See \cite{N3, \S 1}.
We recall that these forms are defined by extension
of the form of the lattice $S$ to $S^\ast$. The discriminant quadratic
form $q_S$ takes values in $\q/2\z$. And $b_S$ in
$\q/\z$.

These invariants $(r(\phi ),a(\phi ),\delta (\phi ))$
define the involution $\phi $ uniquely up to automorphisms of the
unimodular lattice $L$
(see \cite{N3, \S  3}). Also, they define the topology of the fixed part
of the action of $\phi $ on $X(\cc )$ (see
\cite{N2, \S  4} for the
 holomorphic case and \cite{H, appendix} and \cite{N3, \S  3}
 for the antiholomorphic one):
$$
X(\cc )^\phi =\cases
\emptyset ,&\text
{if  $(r(\phi ),a(\phi ),\delta (\phi ))=(10,10,0);$}\\
2T_1, &\text{if $(r(\phi ),a(\phi ),\delta (\phi ))=(10,8,0);$}\\
T_{g(\phi )}\amalg k(\phi)T_0,  &\text{where
$g(\phi )=(22-r(\phi )-a(\phi ))/2,$}\\
&\text{$k(\phi )=(r(\phi )-a(\phi ))/2$ otherwise.}
\endcases
\tag3--3--4
$$
Here $T_g$ is a real orientable compact surface of  genus $g$.

We mention the basic formulae necessary to get (3--3--4).

The invariant $\delta(\phi)$ has the following geometrical sense:
$$
X(\cc )^\phi \sim 0\mod 2\ in\ H_2(X(\cc );\z)\ iff\  \delta (\phi )=0.
\tag3--3--5
$$
In particular, it follows that the invariant $\delta (\phi )=0$ if
$X(\cc )^\phi =\emptyset $.

Using (1--0), from the Lefschetz fixed point formula
(see \cite{Sp, Ch. 5, Sect. 7, Theorem 6}),
$$
r(\phi)=\chi (X(\cc )^\phi )/2+10.
\tag3--3--6
$$
 From \cite{H, Lemma 3.7} (it is the consequence of Smith exact
sequence applied to the involution $\phi $, see \cite{Br, Ch.III}),
$$
\dim H^\ast (X(\cc )^\phi ; \z /2)=
\dim H^\ast (X(\cc ); \z /2)-2a(\phi ),
$$
if $X(\cc )^\phi \not= \emptyset $.
It follows therefore that
$$a(\phi )=
\cases
12-\dim H^\ast (X(\cc )^\phi ;\z /2)/2,
&\text{if $X(\cc )^\phi \not= \emptyset ,$}\\
10,  &\text{if $X(\cc )^\phi =\emptyset .$}
\endcases
\tag3--3--7
$$

For any connected closed surface $F$
(orientable or not) we have the
formula
$$
\chi (F)/4+\dim H^\ast (F;\z /2)/4=1.
$$
Then, from (3--3--6) and (3--3--7), for the number $s(\phi )$ of
connected components of $X(\cc )^\phi $, we get the
following formula:
$$
s(\phi )=
\cases
(r(\phi )-a(\phi ))/2+1 &\text{if $s(\phi )>0$,}\\
(r(\phi )-a(\phi ))/2 &\text{if $s(\phi )=0$.}
\endcases
\tag3--3--8
$$

Since the involution $\tau $ has no fixed points, we get
$$
(r(\tau ),a(\tau ),\delta (\tau))=(10,10,0).
$$
Further, we want to connect the invariants
$$
(r(\sigma),a(\sigma),\delta (\sigma)),\ and\
(r(\tau\sigma),a(\tau\sigma),\delta (\tau\sigma))
$$
of the involutions
$\sigma $ and $\tau \sigma $.

\smallpagebreak

Let $L$ be an even unimodular lattice, $S\subset L$ a primitive
sublattice of $L$ and $\theta \mid S$ an involution of $S$.
In \cite{N4} all genus invariants of an extension of $\theta$
to an involution
$\sigma $ of the unimodular lattice $L$ were
found and studied.
In \cite{N4} it was formulated as a studying of triplets
$$
(L, S, \sigma ).
$$
Here $L$ is an even unimodular lattice, $S$ is a
primitive sublattice of
$L$
(i.e. $L/S$ is free)
and $\sigma $ an involution of $L$ such that $\sigma (S)=S$ and
$\sigma \mid S=\theta$ where $\theta$ is a fixed involution of $S$.
We apply this theory to the triplet $(L,L^\tau , \sigma )$.

We give some invariants from \cite{N4} of the
triplet $(L,L^\tau , \sigma )$ which are necessary for us to prove
our results here. The full investigation of the
invariants from \cite{N4} for
real Enriques surfaces is in \cite{N6}.

\smallpagebreak

At first, we have to study invariants of the action
$\theta =\sigma \mid L^\tau$.

The lattice $L^\tau $ is an even $2$-elementary lattice with invariants
$r(\tau )=\rk L^\tau =10,
\ a(\tau )=\dim (L^\tau )^\ast /L^\tau =10$ and an
even discriminant 2-elementary form $q_{L^\tau }$ (see (3--3--3)).
It follows that the
lattice $L^\tau (1/2)$ is an even unimodular lattice of signature
(1,9). Here we use the following notation:

\definition{Notation 3.3.1} For a lattice $M$, the lattice $M(a)$ is defined
by multiplying the form of the lattice $M$
by $a\in {\bold  Q}$.
\enddefinition

Thus, the restriction of the involution $\sigma $ on
the lattice $L^\tau $ is defined by the action of $\sigma $ on
the unimodular lattice $L^\tau (1/2)$. By (1--6),  we have the
canonical isomorphism of lattices:
$$
\pi^\ast :\  H^2(Y(\cc );\z )/\Tor \cong
L^\tau (1/2).
$$
because $\pi^\ast (x)\cdot \pi^\ast(y)=2x\cdot y$ for
$x,y\in H^2(Y(\cc );\z)$.
The map $\pi^\ast $ is equivariant with respect to the action of $\theta $
and its lifting $\sigma $. Hence
it follows that the action of $\sigma $ on $L^\tau $ is
the same as the action of $\theta$ on
$H^2(Y(\cc );\z)/ \Tor$.
Thus, our notation
$\theta =\sigma \mid L^\tau $
agrees with our previous notation
$\theta $ for an antiholomorphic involution of $Y$.
For the involution  $\theta =\sigma \mid L^\tau$ on the even
unimodular lattice
$L^\tau (1/2)$, we have a similar
triplet of invariants as above (for $\sigma $ and $\tau\sigma $ on $L$):
$$
(r(\theta ), a(\theta ), \delta (\theta )),
\tag3--3--9
$$
where
$$
r(\theta )=\rk  L^{\tau ,\sigma }
$$
(we denote $L^{\tau ,\sigma }=(L^\tau)^\sigma$ and
$L_{\tau ,\sigma}=(L_\tau)_\sigma$);
$$
(L^{\tau ,\sigma }(1/2))^\ast /
L^{\tau ,\sigma }(1/2) \cong (\z /2)^{a(\theta )};
$$
$\delta (\theta ) $ is zero if the discriminant quadratic form
$q_{L^{\tau ,\sigma }(1/2)}$ takes values in $\z /2\z $.
Otherwise, $\delta (\theta )=1.$

Since the lattices $L^\tau ,L^\sigma $ and $L^{\tau \sigma }$ are
hyperbolic, the lattice $L^{\tau ,\sigma }$ is negative and the
lattice $L^\tau _\sigma $ is hyperbolic of
$\rk L^\tau _\sigma =10-r(\theta ).$
Using results of \cite{N3} (or \cite{N4})
one can describe very easily all  possibilities
for triplets of invariants $(r(\theta ), a(\theta ), \delta (\theta ))$.
They are:
$$
\split
(r(\theta ), a(\theta ), \delta (\theta ))=
& (0,0,0),\ (1,1,1),\ (2,2,1),\ (3,3,1),\ (4,2,0),\ (4,4,1),\\
& (5,3,1),\ (5,5,1),\ (6,2,1),\ (6,4,1),\ (7,1,1),\ (7,3,1),\\
& (8,0,0),\ (8,2,0),\ (8,2,1),\ (9,1,1).
\endsplit
\tag3--3--10
$$

We need another interpretation of the invariants
$(r(\theta),a(\theta ),\delta (\theta ))$ using the lattices
$L^{\tau ,\sigma }$ and $L^\tau _\sigma $. We give it below.

For the invariant $r(\theta)$ we evidently have
$r(\theta )=\rk L^{\tau ,\sigma }$, $\rk L^\tau _\sigma =10-r(\theta)$.

We have the canonical
subgroup
$$
\Gamma (L^\tau (1/2))\subset
A_{L^{\tau ,\sigma }(1/2)}\oplus A_{L^\tau _\sigma (1/2)},
$$
where
$\Gamma (L^\tau (1/2))=
L^\tau (1/2)/(L^{\tau ,\sigma }(1/2)\oplus L^\tau _\sigma (1/2)).$
For the discriminant form
$q_{L^{\tau ,\sigma }(1/2)}\oplus q_{L^\tau _\sigma (1/2)}$ on
the group
$A_{L^{\tau ,\sigma }(1/2)}\oplus A_{L^\tau _\sigma (1/2)}$
the subgroup
$\Gamma (L^\tau (1/2))$
is isotropic and
$$
(\Gamma (L^\tau (1/2)))^\perp =\Gamma (L^\tau (1/2))
$$
since the lattice $L^\tau (1/2)$ is unimodular. Moreover,
$$
\Gamma (L^\tau (1/2))\cap A_{L^{\tau ,\sigma }(1/2)}=
\Gamma (L^\tau (1/2))\cap A_{L^\tau _\sigma (1/2)}=\{ 0\} .
$$
Let $p_1$ and $p_2$ be the projections
on
$A_{L^{\tau ,\sigma }(1/2)}$ and $A_{L^\tau _\sigma (1/2)}$
respectively. Then the map
$(p_2 \mid \Gamma (L^\tau (1/2)))(p_1)^{-1}$
is an isomorphism of the discriminant quadratic forms
$q_{L^{\tau ,\sigma }(1/2)}$ and $-q_{L^\tau _\sigma (1/2)}$
on the groups
$A_{L^{\tau ,\sigma }(1/2)}$ and
$ A_{L^\tau _\sigma (1/2)}$
respectively.
One sees very easily by considering an orthogonal
decomposition of 2-adic lattices
$L^{\tau ,\sigma }\otimes \z_2,\ L^\tau _\sigma\otimes \z_2$
as a sum of elementary lattices of rank one or two, that
we have the following identifications where for a group $A$ we
denote $A^{(1)}=Ann (2\ in\ A)$:
$$
A_{L^{\tau ,\sigma }(1/2)}=
(L^{\tau ,\sigma }(1/2))^\ast /L^{\tau ,\sigma }(1/2)=
2(L^{\tau ,\sigma })^\ast /L^{\tau ,\sigma }=
$$
$$
=2A_{L^{\tau ,\sigma }}=\Ker (b_{L^{\tau ,\sigma }}\mid
A_{L^{\tau ,\sigma }}^{\ (1)}).
$$
The same is true for the lattice
$L^\tau _\sigma $. We have
$$
A_{L^\tau _\sigma (1/2)}=
(L^\tau _\sigma (1/2))^\ast /L^\tau _\sigma (1/2)=
2(L^\tau _\sigma )^\ast /L^\tau _\sigma =
$$
$$
=2A_{L^\tau _\sigma }=\Ker (b_{L^\tau _\sigma }\mid
A_{L^\tau _\sigma }^{\ (1)}).
$$
Let us denote
$$
\Gamma (\sigma )_+=
2A_{L^{\tau ,\sigma }}=\Ker (b_{L^{\tau ,\sigma }}
\mid A_{L^{\tau ,\sigma }}^{\ (1)})
$$
and
$$
\Gamma (\sigma )_-=
2A_{L^\tau _\sigma }=\Ker (b_{L^\tau _\sigma }\mid
A_{L^\tau _\sigma }^{\ (1)}).
$$
Moreover, arguing similarly, we see that
$$
q_{L^{\tau, \sigma }}\mid \Gamma (\sigma )_+=0\ \  iff\ \
\delta (\theta )=0.
$$
And
$$
q_{L^\tau _\sigma }\mid \Gamma (\sigma )_-=0\ \ iff\ \
\delta (\theta )=0.
$$
Thus, we get the following interpretation:

\proclaim
{Lemma 3.3.1} Let $Y$ be an arbitrary real Enriques surface  with
an antiholomorphic involution $\theta$. Let
$(r(\theta ),a(\theta ),\delta (\theta ))$
are invariants described above of the action of $\theta$ on
the lattice $H^2(Y(\cc );\z )/\Tor$
(equivalently, invariants of the action of the
lifting $\sigma $ of $\theta $ on the
lattice $\L^\tau (1/2)$).

Then:
$r(\theta)=\rk L^{\tau ,\sigma }$ and $\rk L^\tau _\sigma =10-r(\theta)$.

For the subgroup
$$
\Gamma _{12}=L^\tau /(L^{\tau ,\sigma }\oplus L^\tau _\sigma )
\subset
A_{L^{\tau ,\sigma }}\oplus A_{L^\tau _\sigma }
$$
and the projections
$p_1$ and $p_2$ on the groups
$A_{L^{\tau ,\sigma }}$ and $A_{L^\tau _\sigma }$ respectively we have:
$$
p_1:(\Gamma _{12})\cong \Gamma (\sigma )_+=\Ker
b_{L^{\tau ,\sigma }}\mid
A_{L^{\tau ,\sigma }}^{\ (1)}
\cong (\z /2)^{a(\theta )}
$$
and
$$
p_2:(\Gamma _{12})\cong \Gamma (\sigma )_-=\Ker
b_{L^\tau _\sigma }\mid
A_{L^\tau _\sigma }^{\ (1)}
\cong (\z /2)^{a(\theta )}
$$
(where $A^{(1)}=Ann(2)$ for an abelian group $A$).

The invariant $\delta (\theta )=0$ iff
$q_{L^{\tau ,\sigma }}\mid
\Gamma(\sigma )_+ = 0$,
equivalently,
$q_{L^\tau _\sigma }\mid
\Gamma (\sigma)_- = 0).$

Besides, we have the following property:
$$
b_{L^{\tau ,\sigma }}(x,x)=0\ for\ any\ x\in
A_{L^{\tau ,\sigma }}^{\ (1)};
$$
and
$$
b_{L^\tau _\sigma }(x,x)=0\ for\ any\ x\in
A_{L^\tau _\sigma }^{\ (1)}.
$$
\endproclaim

\demo{Proof} We should only explain the last statement.
This is true because the lattice
$L^\tau (1/2)$ is even:
$z^2\equiv 0 \mod 2$ for any $z\in L^\tau (1/2)$.
\enddemo

Thus, we studied the action $\theta =\sigma \mid L^\tau =
\tau\sigma \mid L^\tau$ introducing
the invariants \linebreak
$(r(\theta),a(\theta),\delta(\theta))$ above.
Now we fix these invariants (one of triplets
(3--3--10)) and want to give some invariants
of extension $\sigma $ of $\theta $
on the unimodular lattice $L$. Thus, we
choose an involution between $\sigma $ and $\tau\sigma$.

We consider the decomposition
$L\supset L^\sigma \oplus L_\sigma $ of finite index.
Up to automorphisms of $L$, this is defined by the invariants

$$
(r(\sigma), a(\sigma),\delta(\sigma))
\tag3--3--11
$$
we had introduced above (see (3--3--1), (3--3--2), (3--3--3)).
These invariants are the most important invariants of the
extension $\sigma $ of $\theta$. But there are some other invariants we want
to describe.

Let us consider the corresponding
decomposition of the discriminant groups
$A_{L^\sigma }\oplus A_{L_\sigma }$
with the discriminant form
$q_{L^\sigma }\oplus q_{L_\sigma }.$ The subgroup
$\Gamma (L^\sigma , L_\sigma )=L/(L^\sigma \oplus L_\sigma )\subset
A_{L^\sigma }\oplus A_{L_\sigma }$ is isotropic,
$\Gamma (L^\sigma , L_\sigma )^\perp =
\Gamma (L^\sigma , L_\sigma )$
and
$\Gamma (L^\sigma , L_\sigma )
\cap A_{L^\sigma }=\Gamma (L^\sigma , L_\sigma )\cap
A_{L_\sigma }=\{ 0\} $. It follows that
$\Gamma (L^\sigma , L_\sigma )$ is the graph of the isomorphism
$\gamma (L^\sigma , L_\sigma ):
q_{L^\sigma }\cong -q_{L_\sigma }$
of the discriminant quadratic forms. Using this isomorphism, we
identify discriminant groups and quadratic forms:
$$
q(\sigma )=q_{L^\sigma }=-q_{L_\sigma },\
A_{q(\sigma )}=A_{L^\sigma }=A_{L_\sigma }.
\tag3--3--12
$$

We now consider the following
decomposition of lattices,
$L^{\tau, \sigma }\oplus L_\sigma $,
and the corresponding decomposition of
the discriminant groups
$A_{L^{\tau, \sigma }}\oplus A_{L_\sigma }$
with the  discriminant form
$q_{L^{\tau, \sigma }}\oplus q_{L_\sigma }$. The subgroup
$\Gamma (L^{\tau, \sigma },L_\sigma )=
(L^\sigma _\tau )^\perp /(L^{\tau, \sigma }\oplus L_\sigma )\subset
A_{L^{\tau, \sigma }}\oplus A_{L_\sigma }$
is isotropic,
$\Gamma (L^{\tau, \sigma },L_\sigma )\cap A_{L^{\tau, \sigma }}=
\Gamma (L^{\tau, \sigma },L_\sigma )\cap A_{L_\sigma }=\{ 0\} .$
Thus, the subgroup
$\Gamma (L^{\tau, \sigma },L_\sigma )$
is the graph of the embedding
$\gamma (L^{\tau, \sigma },L_\sigma )$
of a subgroup
$H(\sigma )_+ \subset A^{\ (1)}_{L^{\tau, \sigma }}$ into the group
$A_{q(\sigma )}=A_{L_\sigma }$
with the discriminant form $q(\sigma )$:
$$
H(\sigma )_+ \subset A^{\ (1)}_{L^{\tau, \sigma }},\ \
\gamma (L^{\tau, \sigma },L_\sigma ):H(\sigma)_+\hookrightarrow
A_{q(\sigma )}.
\tag3--3--13
$$
Evidently,
$\gamma (L^{\tau, \sigma },L_\sigma )$
gives the embedding of quadratic forms
$$
\gamma (L^{\tau, \sigma },L_\sigma ):
q_{L^{\tau, \sigma }}\mid H(\sigma )_+ \to q(\sigma )\mid A_{q(\sigma)}.
\tag3--3--14
$$
We repeat the same construction for
$L^\tau _\sigma \oplus L^\sigma $.
We have the decomposition
$A_{L^\tau _\sigma }\oplus A_{L^\sigma }$
with the discriminant form
$q_{L^\tau _\sigma }\oplus q_{L^\sigma }$.
The subgroup
$\Gamma (L^\tau _\sigma , L^\sigma )=
(L_{\tau ,\sigma })^\perp /(L^\tau _\sigma \oplus L^\sigma )
\subset
A_{L^\tau _\sigma }\oplus A_{L^\sigma }$
is isotropic,
$\Gamma (L^\tau _\sigma ,L^\sigma )\cap A_{L^\tau _\sigma }=
\Gamma (L^\tau _\sigma ,L^\sigma )\cap A_{L^\sigma }=\{ 0\} .$
The subgroup
$\Gamma (L^\tau _\sigma ,L^\sigma )\subset
A_{L^\tau _\sigma }\oplus A_{L^\sigma }$
is the graph of the embedding
$\gamma (L^\tau _\sigma ,L^\sigma )$
of a subgroup
$H(\sigma )_- \subset
A^{\ (1)}_{L^\tau _\sigma }$ into the group
$A_{q(\sigma )}=A_{L^\sigma }$
with the discriminant form $q(\sigma )$:
$$
H(\sigma )_- \subset A^{\ (1)}_{L^\tau _\sigma },\ \
\gamma (L^\tau _\sigma ,L^\sigma ):
H(\sigma )_- \hookrightarrow A_{q(\sigma )}.
\tag3--3--15
$$
Evidently,
$\gamma (L^\tau _\sigma ,L^\sigma )$
gives the embedding of quadratic forms
$$
\gamma (L^\tau _\sigma ,L^\sigma ):
-q_{L^\tau _\sigma }\mid H(\sigma )_- \to q(\sigma )\mid A_{q(\sigma)}.
\tag3--3--16
$$
Above, we had defined subgroups
$\Gamma (\sigma )_+\subset A^{\ (1)}_{L^{\tau ,\sigma }}$
and
$\Gamma (\sigma )_-\subset A^{\ (1)}_{L^\tau _\sigma }.$
One can see very easily (see \cite{N4, \S  1}) that
$\Gamma (\sigma )_\pm \subset H(\sigma )_\pm$ and
$$
\split
& \gamma (L^{\tau ,\sigma },L_\sigma )
(H(\sigma )_+)\cap
\gamma (L^\tau _\sigma ,L^\sigma )(H(\sigma )_-)=\\
& = \gamma (L^{\tau ,\sigma },L_\sigma )(\Gamma (\sigma )_+)
= \gamma (L^\tau _\sigma ,L^\sigma )(\Gamma (\sigma )_-).
\endsplit
\tag3--3--17
$$

The triplet of the basic invariants
$(r(\sigma ),a(\sigma ),\delta(\sigma ))$ and
subgroups $H(\sigma )_+ \subset
A^{\ (1)}_{L^{\tau, \sigma }}$ and
$H(\sigma )_- \subset
A^{\ (1)}_{L^\tau _\sigma }$ with
the pairing between them defined by
embeddings $\gamma (L^{\tau ,\sigma },L_\sigma )$ and
$\gamma (L^\tau _\sigma ,L^\sigma )$  respectively
and the discriminant form
$q(\sigma )$ are the most important invariants from \cite{N4}
of the triplet $(L, L^\tau , \sigma  )$.

Other more delicate invariants of the triplet $(L, L^\tau ,\sigma )$
can be found in  \cite{N4}.
We will use some of them in Sect. 3.4.
A complete description of invariants from \cite{N4} for
real Enriques surfaces can be found in \cite{N6}.

\definition{Notation 3.3.2} To simplify our notations, we will identify
groups $H(\sigma )_+$ and $H(\sigma )_-$
with their images
$\gamma (L^{\tau, \sigma },L_\sigma )(H(\sigma)_+)$ and
$\gamma (L^\tau _\sigma ,L^\sigma )(H(\sigma )_-)$ in $A_{q(\sigma )}$.
See (3--3--13) and (3--3--15). Thus, by (3--3--17), we have
$$
H(\sigma )_+\cap H(\sigma )_-=\Gamma (\sigma )_\pm
$$
(we use similar identification for groups $\Gamma (\sigma )_\pm$).
Using the discriminant quadratic form $q(\sigma )$ on
$A_{q(\sigma )}$, we then have a bilinear pairing between $H(\sigma )_+$ and
$H(\sigma )_-$ and can consider the orthogonal complements
$H(\sigma )_\pm^\perp$ to $H(\sigma )_\pm$ in $A_{q(\sigma )}$.
We remark that by (3--3--14) and (3--3--16), we have
$$
q_{L^{\tau, \sigma }}\mid H(\sigma )_+=q(\sigma )\mid H(\sigma )_+,\
-q_{L^\tau _\sigma}\mid H(\sigma )_-=q(\sigma )\mid H(\sigma )_-.
$$
\enddefinition

The following very important
and non-trivial Lemma connects the invariants \linebreak
$(r(\sigma ),a(\sigma ),\delta (\sigma ))$ and
$(r(\tau \sigma ),a(\tau \sigma ),\delta (\tau \sigma ))$ of the
involutions $\sigma $ and $\tau \sigma $. We will not use the
third statement of the Lemma here, and include it for the
sake of completeness.

\proclaim
{Lemma 3.3.2} Let $Y$ be a real Enriques surface such that there are two
liftings $\sigma $ and $\tau \sigma $ of the antiholomorphic involution
$\theta $ of $Y$ to the antiholomorphic involutions of the $K3$-surface
$X(\cc )$ (it is true for example if
$Y(\r )\not= \emptyset $).

Then
$$
r(\sigma )+r(\tau \sigma )=12+2r(\theta );
$$
$$
\split
&a(\tau \sigma )-a(\sigma )=\\
&=10+2a(\theta )-2\dim H(\sigma )_+
-2\dim (H(\sigma )_+)^\perp
\cap H(\sigma )_-;
\endsplit
$$
$$
\delta (\sigma )+\delta (\tau \sigma )\equiv \delta (\theta ) \mod 2.
$$
\endproclaim

\demo{Proof} We consider the decomposition
$$
L^{\tau ,\sigma }\oplus L^\tau _\sigma \oplus
L_\tau ^\sigma \oplus L_{\tau, \sigma}\subset L
$$
of finite index
and the corresponding direct sum of the discriminant groups
$$
A_{L^{\tau ,\sigma }}\oplus A_{L^\tau _\sigma }\oplus
A_{L_\tau ^\sigma }\oplus A_{L_{\tau, \sigma}}
$$
with the discriminant quadratic form
$$
q_{L^{\tau ,\sigma }}\oplus q_{L^\tau _\sigma }\oplus
q_{L_\tau ^\sigma }\oplus q_{L_{\tau, \sigma}}
$$
and isotropic subgroup
$$
\Gamma =L/(L^{\tau ,\sigma }\oplus L^\tau _\sigma \oplus
L_\tau ^\sigma \oplus L_{\tau, \sigma})\subset
A_{L^{\tau ,\sigma }}\oplus A_{L^\tau _\sigma }\oplus
A_{L_\tau ^\sigma }\oplus A_{L_{\tau, \sigma}}
$$
Let $p_1, p_2, p_3, p_4$ be the corresponding projections on the
$A_1=A_{L^{\tau ,\sigma }}$,  $A_2=A_{L^\tau _\sigma }$,
$A_3=A_{L_\tau ^\sigma }$, $A_4=A_{L_{\tau, \sigma}}$
respectively. We denote by $q_i$
respectively the corresponding
discriminant form on the group $A_i$. Evidently,
$\Gamma \cap A_i=\{ 0\} , 1\le i \le 4.$

Let $\Gamma _{ij}=\Gamma \cap (A_i\oplus A_j)$
for $1\le i<j\le 4$. We use a similar definition for
$\Gamma _{ijk}$, where $1\le i<j<k\le 4$. We remark that
Lemma 3.3.1 gives the description of the subgroup
$\Gamma _{12}$.

Evidently, $\rk L^\sigma =
\rk L^{\tau, \sigma }+\rk L_\tau ^\sigma $
and
$\rk L^{\tau \sigma }=
\rk L^{\tau, \sigma }+\rk L_{\tau ,\sigma }.$
Thus, it follows that
$\rk L^\sigma+\rk L^{\tau \sigma }=
2\rk L^{\tau, \sigma }+\rk L_\tau ^\sigma +\rk L_{\tau ,\sigma }=
2\rk L^{\tau, \sigma }+\rk L_\tau =2r(\theta )+12.$
This proves the first
statement.

We now prove the second part. The following statement will be useful to us.

\proclaim
{ Proposition 3.3.3} Let
$S_1\oplus S_2 \oplus S_3\subset L$
be an orthogonal decomposition up to finite index of an unimodular
lattice $L$ with all sublattices $S_i, 1\le i\le 3$
are primitive in $L$. Let
$\Gamma (L)=L/(S_1\oplus S_2\oplus S_3)\subset
A_{S_1}\oplus A_{S_2}\oplus A_{S_3}$ be the corresponding isotropic
subgroup with respect to the discriminant form
$q_{S_1}\oplus q_{S_2}\oplus q_{S_3}$. Let
$\Gamma (S_i,S_j)=\Gamma (L)\cap (A_{S_i}\oplus A_{S_j})$ where
$1\le i<j\le 3.$

Then the orthogonal complement to $\Gamma (S_1,S_2)$ in
$\Gamma (L)$ with respect to the form $q_{S_2}$ is equal to
$\Gamma (S_2,S_3)+\Gamma (S_1,S_3)$.
\endproclaim

\demo{Proof} Evidently, these subgroups are orthogonal to one another.
Let $p_i$ be the projection on $A_{S_i}$.

We prove that
$p_2(\Gamma (S_1,S_2))$ and $p_2(\Gamma (S_2,S_3))$ are the orthogonal
complements to one another with respect to the form $q_{S_2}$.
The projection $p_2$ is an embedding of both these subgroups since
the sublattices $S_i$ are primitive. It follows that
$$\sharp A_{S_1}= \sharp A_{(S_1)^\perp _L}=
(\sharp A_{S_2}\sharp A_{S_3})/(\sharp \Gamma (S_2,S_3))^2$$
and
$$\sharp A_{S_3} =\sharp A_{(S_3)^\perp _L}=
(\sharp A_{S_1}\sharp A_{S_2})/(\sharp \Gamma (S_1,S_2))^2.$$
It follows therefore that
$$\sharp A_{S_2}=\sharp \Gamma (S_1,S_2)\sharp \Gamma (S_2,S_3).$$
This equality is equivalent to the statement we had to prove
since the bilinear form of the form $q_{S_2}$ is
non-degenerate.

We now consider the projection $p_2:\Gamma (L)\to A_{S_2}$. Evidently,
$ker\ p_2=\Gamma (S_1,S_3)$. From the result  proved above,
$(p_2)^{-1}(p_2(\Gamma (S_1,S_2))^\perp )=
\Gamma (S_2,S_3)+\Gamma (S_1,S_3)$, where we take the orthogonal complement
in $A_{S_2}$ using the bilinear form of $q_{S_2}$.
Hence the Proposition follows.

\enddemo

 From the definition of $r({\theta })$ and $a({\theta })$,
$$\sharp A_{L^{\tau ,\sigma }}=
4^{a(\theta )}2^{r(\theta )-a(\theta )}=
2^{r(\theta )+a(\theta )};
\tag3--3--18
$$
$$\sharp A_{L^\tau _\sigma }=
4^{a(\theta )}2^{10-r(\theta )-a(\theta )}
=2^{10-r(\theta )+a(\theta )};
\tag3--3--19
$$
$$
\split
\sharp A_{L_\tau ^\sigma }=\sharp A_{(L_\tau ^\sigma )^\perp _L}&=
(\sharp A_{L^{\tau ,\sigma }}\sharp A_{L_\sigma })/(\sharp H(\sigma )_+)^2\\
&=2^{r(\theta )+a(\theta )+a(\sigma )-2\dim H(\sigma )_+};
\endsplit
\tag3--3--20
$$
$$
\split
\sharp A_{L_{\tau ,\sigma }}=\sharp A_{(L_{\tau ,\sigma })^\perp_L}&=
(\sharp A_{L^\tau _\sigma }\sharp A_{L^\sigma })/(\sharp H(\sigma )_-)^2\\
&=2^{10-r(\theta )+a(\theta )+a(\sigma )-2\dim H(\sigma )_-}.
\endsplit
\tag3--3--21
$$

By the projection $p_1\oplus p_2$, we have the identifications
$$q(\sigma )\mid A_{q(\sigma )}
=q_1\oplus q_3\mid \Gamma /(\Gamma _{13}+\Gamma _{24});$$
$$\gamma (L^{\tau ,\sigma },L_\sigma )H(\sigma )_+
=\Gamma _{124}/\Gamma _{24}\subset
\Gamma /(\Gamma _{13}+\Gamma _{24});$$
$$\gamma (L^\tau _\sigma ,L^\sigma )H(\sigma )_-
=\Gamma _{123}/\Gamma _{13}\subset
\Gamma /(\Gamma _{13}+\Gamma _{24}).$$
By Proposition 3.3.3 and this identification,
$$(\gamma (L^{\tau ,\sigma },L_\sigma )H(\sigma )_+)^\perp _{q(\sigma )}=
(\Gamma _{13}+\Gamma _{234})/(\Gamma _{13}+\Gamma _{24}).$$
We hence have
$$(\gamma (L^{\tau ,\sigma },L_\sigma )H(\sigma )_+)^\perp _{q(\sigma )}
\cap \gamma (L^\tau _\sigma ,L^\sigma )H(\sigma )_-=$$
$$=(\Gamma _{13}+\Gamma _{234})\cap (\Gamma _{24}+\Gamma _{123})
/(\Gamma _{13}+\Gamma _{24})=$$
$$=(\Gamma _{13}+\Gamma _{24}+
\Gamma _{234}\cap \Gamma _{123})/(\Gamma _{13}+\Gamma _{24})\cong
\Gamma _{23}.$$

Thus,
$$\sharp (\gamma (L^{\tau ,\sigma },L_\sigma )H(\sigma )_+)^\perp _{q(\sigma )}
\cap \gamma (L^\tau _\sigma ,L^\sigma )H(\sigma )_-=
\sharp \Gamma _{23}.$$
On the other hand, by definition, and calculations above

$$2^{a(\tau \sigma )}=\sharp A_{L^{\tau \sigma }}=
\sharp A_{L_{\tau \sigma }}=
\sharp A_{L^\tau _\sigma }\sharp A_{L^\sigma _\tau }/(\sharp \Gamma _{23})^2=$$
$$=2^{10-r(\theta )+a(\theta )}
2^{r(\theta )+a(\theta )+a(\sigma )-2\dim H(\sigma )_+}2^{-2\dim  \Gamma
_{23}}.$$
This gives the proof of the second statement of the  Lemma.

We now prove the last statement. As before,
$q(\sigma )\mid A_{q(\sigma )}=
q_1\oplus q_3\mid \Gamma /(\Gamma _{13}+\Gamma _{24}).$ It follows that
$$\delta (\sigma )=0\ iff\  q_1+q_3\equiv 0\mod 1.$$
(This must be true on the group $\Gamma $.) The same holds for the involution
$\tau\sigma$:
$$\delta (\tau\sigma )=0\ iff\ q_1+q_4\equiv 0\mod 1.$$
We remark that
$$q_1+q_2+q_3+q_4=0$$
since the subgroup $\Gamma $ is isotropic. Moreover, we know that
$\delta (\tau )=0$. If follows that
$$q_1+q_2\equiv q_3+q_4\equiv 0 \mod 1.$$
It follows very easily that if $\delta (\sigma )=0$ then
$\delta (\tau \sigma )=0$ iff $2q_1\equiv 0 \mod 1$. On the other
hand, one can see very easily that $2q_1\equiv 0 \mod 1$ iff
$\delta (\theta )=0$. By symmetry, if $\delta (\tau\sigma )=0$,
then $\delta (\sigma )=0$ iff $\delta (\theta )=0$. This proves the
last statement.
\enddemo

 From Lemma 3.2.1, formula (3--3--18) and
Lemma 3.3.2, we get the following very important

\proclaim{Theorem 3.3.4}
Let $Y$ be a real Enriques surface such that there are two
liftings $\sigma $ and $\tau \sigma $ of the antiholomorphic involution
$\theta $ of $Y$ to antiholomorphic involutions of the $K3$-surface
$X(\cc )$ (this is true for example if
$Y(\r )\not= \emptyset $).
Let $s(\sigma )$ and  $s(\tau \sigma )$ be the number of connected
components of the real parts $X_\sigma(\r ) $ and
$X_{\tau\sigma}(\r )$ respectively.

Then
$$
\split
s(\sigma )+s(\tau \sigma )
& =2+(r(\sigma )+r(\tau \sigma ))/2-(a(\sigma )+a(\tau \sigma ))/2\\
& =3+r(\theta )-a(\theta )-a(\sigma )+\dim H(\sigma )_+
+ \dim (H(\sigma )_+)^\perp
\cap H(\sigma )_-
\endsplit
$$
if both $s(\sigma )>0$ and
 $s(\tau \sigma )>0.$

If either $s(\sigma )=0$ and $s(\tau \sigma )>0$
or $s(\sigma )>0$ and $s(\tau \sigma )=0$ then
$$
\split
s(\sigma )+s(\tau \sigma )
& = 1+(r(\sigma )+r(\tau \sigma ))/2-(a(\sigma )+a(\tau \sigma ))/2\\
& =2+r(\theta )-a(\theta )-a(\sigma )
 + \dim  H(\sigma )_+ +
\dim (H(\sigma )_+)^\perp
\cap H(\sigma )_-.
\endsplit
$$

If $s(\sigma )=s(\tau \sigma )=0$, then
$$
\split
0=s(\sigma )+s(\tau \sigma )
& = (r(\sigma )+r(\tau \sigma ))/2-(a(\sigma )+a(\tau \sigma ))/2\\
&=1+r(\theta )-a(\theta )-a(\sigma )+ \dim  H(\sigma )_+ +
\dim (H(\sigma )_+)^\perp
\cap H(\sigma )_-.
\endsplit
$$

Besides, $s(\sigma )+s(\tau \sigma )=s_{nor}+2s_{or}$ where
$s_{or}$ and $s_{nor}$ is the number of orientable and non-orientable
 connected components of $Y(\r )$ respectively.
\endproclaim

\subhead
3.4. Calculation of group cohomology invariants
\endsubhead

Here we calculate
group cohomology using the invariants above.

First of all, we have the following simple

\proclaim{Proposition 3.4.1} Let $Y$ be an
arbitrary real Enriques surface . Then
$$
\dim (\Pic Y_\cc )^\theta /2(\Pic Y_\cc )^\theta =
11-r(\theta ).
$$

\endproclaim

\demo{Proof} See the proof of the formula (2--2) in \S 2.
\enddemo

We use the following
\proclaim
{Proposition 3.4.2} Let $Y$ be a real Enriques surface and suppose
there are two
liftings $\sigma $ and $\tau \sigma $ of the antiholomorphic involution
$\theta $ of $Y$ to antiholomorphic involutions of the $K3$-surface
$X(\cc )$ (for example it is true if
$Y(\r )\not= \emptyset $).

 Then there are the canonical
isomorphisms:
$$
\split
H^2(X(\cc );\z /2)^\tau &\cong
(H^2(X(\cc ); \z )^\tau \oplus
H^2(X(\cc ); \z )_\tau )/2H^2(X(\cc ); \z )\\
&\cong
H^2(X(\cc ); \z )_\tau /2H^2(X(\cc ); \z )_{\tau };
\endsplit
$$
$$
H^2(X(\cc );\z /2)^{\tau ,\sigma } \cong
(H^2(X(\cc );\z )_\tau /2H^2(X(\cc );\z )_{\tau })^\sigma \cong
$$
$$
\cong (H^2(X(\cc );\z )_\tau ^\sigma \oplus
H^2(X(\cc );\z )_{\tau ,\sigma })
/2H^2(X(\cc );\z )_{\tau };
$$
and
$$
\split
\dim H^2(X(\cc );\z /2)^{\tau ,\sigma }
&=\dim (H^2(X(\cc );\z )_\tau /2H^2(X(\cc );\z )_{\tau })^\sigma \\
&=12-a(\theta )-a(\sigma )+\dim H(\sigma )_+ +\dim  H(\sigma )_-.
\endsplit
$$
$$
(H^2(X(\cc );\z )^\tau /2H^2(X(\cc ); \z )^{\tau })^\sigma \cong
$$
$$
\cong (H^2(X(\cc );\z )^{\tau ,\sigma }\oplus
H^2(X(\cc );\z )^\tau _\sigma )
/2H^2(X(\cc );\z )^\tau \ ,
$$
and
$$
\dim (H^2(X(\cc );\z )^\tau /2H^2(X(\cc );\z )^\tau )^\sigma =
10-a(\theta ).
$$
\endproclaim

\demo{Proof}
As above, let $L=H^2(X(\cc );\z )$. Then,
since a $K3$-surface has no
torsion in cohomology (see the beginning of \S 1),
$H^2(X(\cc );\z /2)=L/2L$ and
we should calculate $(L/2L)^\tau $. We claim (cl. \cite{Ha1}) that
$$(L/2L)^\tau =(L^\tau \oplus L_\tau )/2L.$$
Let $x=(x_++x_-)/2 \mod 2L \in (L/2L)^\tau $ where
$x_+\in L^\tau ,x_-\in L_\tau $. We have $x-\tau (x)=x_-\in 2L$. It
then follows
that $x=y_++y_-$ where $y_+\in L^\tau ,y_-\in L_\tau $.
Thus, we have proved the
statement.

By (1--0) and since $a(\tau)=10$, we have:
$\dim (L^\tau \oplus L_\tau )/2L=
22-\dim L/(L^\tau \oplus L_\tau ) =12.$ But the group
$(L^\tau \oplus L_\tau )/2L$ contains the subgroup $L_\tau /2L_\tau $
which has dimension 12 too.
Thus, these groups are isomorphic.

The same proof for $L_\tau $ and the involution $\sigma $ of $L_\tau $
gives the second statement of the Proposition.

We now prove the third statement.
We follow the notation in the proof of  Lemma 3.3.2.
We have the sequence of embeddings:
$$
L_\tau \supset  L_\tau ^\sigma \oplus L_{\tau ,\sigma }\supset 2L_\tau
$$
where
$$
\dim H^2(X(\cc );\z /2)^{\tau ,\sigma }=
\dim (L_\tau ^\sigma \oplus L_{\tau ,\sigma })/ 2L_\tau .
$$
It follows, therefore that
$$
2^{\dim (L_\tau ^\sigma \oplus L_{\tau ,\sigma })/ 2L_\tau }=
2^{12}/\sharp (L_\tau /(L_\tau ^\sigma \oplus L_{\tau ,\sigma })).
$$
On the other hand,
$$
\sharp A_{L_\tau }=2^{10}=
(\sharp A_{L_\tau ^\sigma }\sharp A_{L_{\tau, \sigma }})/
\sharp (L_\tau /(L_\tau ^\sigma \oplus L_{\tau ,\sigma }))^2.
$$
 From the calculations above
of
$\sharp A_{L_\tau ^\sigma }$ and $\sharp A_{L_{\tau, \sigma }}$
(see (3--3--20) and (3--3--21)),
 we get that
$$
\dim L_{\tau }/(L_\tau ^\sigma \oplus L_{\tau ,\sigma })
=a(\theta )+a(\sigma )-\dim  H(\sigma )_+-\dim  H(\sigma )_-.
$$
Hence the third statement follows.

The proof of statements about the group
$(H^2(X(\cc );\z )^\tau /2H^2(X(\cc );\z )^{\tau })^\sigma  $
is similar. We remark that from Lemma 3.3.1
$$
\dim L^\tau /(L^{\tau ,\sigma }\oplus L^\tau _\sigma )=a(\theta ).
$$
Besides, $\rk  L^\tau =10$.
\enddemo

We need some information about complex Enriques surfaces.

 From the universal coefficients formula (see \cite{Sp}),
we have a filtration
$$
\Tor H^2(Y(\cc );\z )\subset H^2(Y(\cc );\z )\otimes \z /2
\subset  H^2(Y(\cc );\z /2),
\tag3--4--1
$$
where by (1--4) and (1--5),
$$
\dim \Tor  H^2(Y(\cc );\z )=1,\  \dim H^2(Y(\cc );\z )\otimes \z /2=11,\
\dim H^2(Y(\cc );\z /2)=12.
$$

Using Smith exact sequence (see \cite {Br, Ch. III, Sect. 3}) and
(3--4--1) and (1--0), we have the exact sequence
$$
\split
0&\to
\Tor H^2(Y(\cc );\z )\to H^2(Y(\cc );\z /2)@> \pi ^\ast >>
H^2(X(\cc );\z /2)\\
&@>\pi _\ast >>H^2(Y(\cc );\z )\otimes \z /2\to 0.
\endsplit
\tag3--4--2
$$
This exact sequence is very important for us, and we clarify.
Smith exact sequence for the involution $\tau $ of $X(\cc )$
gives an exact sequence (since $\tau$ has no a fixed point):
$$
\split
0&\to
H^1(Y(\cc );\z/2)
@> \partial^1 >> H^2(Y(\cc );\z /2)@> \pi ^\ast >>
H^2(X(\cc );\z /2)\\
&@>\pi _\ast >>H^2(Y(\cc );\z/2 ) @> \partial^2 >>
H^3(Y(\cc );\z/2)\to 0 .
\endsplit
\tag3--4--2'
$$
(We used here that $H^1(X(\cc );\z/2)=H^3(X(\cc );\z/2)=0$.)
Since a $K3$-surface has no torsion in cohomology, by (3--4--1), we
evidently get that
$\Tor H^2(Y(\cc );\z)$ lies in a kernel of the homomorphism
$\pi^\ast :H^2(Y(\cc );\z /2)\to H^2(X(\cc );\z /2)$.
By (1--4) and (1--5), we have $\dim \Tor H^2(Y(\cc );\z)=
\dim H^1(Y(\cc );\z /2)=1$. Thus, in (3--4--2) and (3--4--2'), the
kernels of $\pi^\ast$ are identified: $\Tor H^2(Y(\cc );\z)=
H^1(Y(\cc );\z/2)\cong \z/2$. Since a $K3$-surface has no torsion in
cohomology, the homomorphism
$\pi _\ast : H^2(X(\cc );\z /2) \to  H^2(Y(\cc );\z/2 )$
is the tensor product by $\z/2$
of the corresponding homomorphism over $\z$. It follows that
$\pi_\ast (H^2(X(\cc );\z /2))\subset
H^2(Y(\cc );\z)\otimes \z/2 \subset H^2(Y(\cc );\z /2)$. Since
$\dim H^3(Y(\cc );\z /2)=1$ and
$\dim H^2(Y(\cc );\z /2)-\dim H^2(Y(\cc );\z)\otimes \z/2 = 1$,
from the (3--4--2'), we get \linebreak
$\pi_\ast (H^2(X(\cc );\z /2))= H^2(Y(\cc );\z)\otimes \z/2$. Thus,
from the exact sequence (3--4--2'), the exact sequence
(3--4--2) follows.

For ${\Cal L}=L/2L$ where (we recall)
$L=H^2(X(\cc );\z)$ and
${\Cal F}=Im\ \pi^\ast $
we have a filtration of subgroups (we use (1--6)):
$$
{\Cal L}^{(\tau )}=L^\tau /2L^\tau \subset {\Cal F} \subset
{\Cal L}^\tau ={\Cal L}_{(\tau )}=L_\tau /2L_\tau  \subset {\Cal L}=L/2L,
\tag3--4--3
$$
where $\dim {\Cal L}^{(\tau )}=10, \dim {\Cal F}=11,
\dim {\Cal L}_{(\tau )}=12,
\dim {\Cal L}=22$.
Then (3--4--2) and (3--4--3) are connected by the canonical isomorphisms
$$
(H^2(Y(\cc );\z )\otimes \z /2)/\Tor H^2(Y(\cc );\z )\cong
\pi ^\ast (H^2(Y(\cc );\z )\otimes \z /2)={\Cal L}^{(\tau )};
\tag3--4--4
$$
$$
H^2(Y(\cc );\z /2)/\Tor H^2(Y(\cc );\z )\cong
\pi ^\ast (H^2(Y(\cc );\z /2))={\Cal F};
\tag3--4--5
$$
$$
H^2(Y(\cc );\z )\otimes \z /2=\pi _\ast ({\Cal L})\cong {\Cal L}/{\Cal F}.
\tag3--4--6
$$

We want to find $\Cal F$. Let $\z _2$ be the ring of 2-adic integers.
The lattice $L_\tau \otimes \z _2\cong
U(1)\oplus 4U(2)$ where
$U=\langle \matrix 0 & 1\\
1 & 0\\
\endmatrix \rangle .
$
It follows that the quadratic form $x^2/2\mod 2$ with symmetric bilinear
form $x\cdot y \mod 2$ ($x,y \in L_\tau )$ on ${\Cal L}_{(\tau )}$ has
10-dimensional kernel ${\Cal L}^{(\tau )}$ and the induced quadratic form
($x^2/2\mod 2$) on
$H={\Cal L}_{(\tau )}/{\Cal L}^{(\tau )}=\{ 0, g_1, g_2, f\} $ has
values $(g_1)^2/2
=(g_2)^2/2=0 \mod 2$ and $f^2/2=1\mod 2$. The subgroup
$F={\Cal F}/{\Cal L}^{\tau }\subset
{\Cal L}_{(\tau )}/{\Cal L}^{(\tau )}=H$
has order two. The following
statement characterizes this subgroup $F$ since
the element $f$ with $f^2/2=1\mod 2 $ is unique.

\proclaim{Lemma 3.4.3}
$F={\Cal F}/{\Cal L}^{\tau }=\{ 0,f\} $, where $f^2/2=1 \mod 2$. In
other words, for $x\in L$ such that $x+2L\in {\Cal F}$, but
$x+2L\notin L^\tau +2L$ we have: $x^2/2=1\mod 2$.
\endproclaim

\demo{Proof} It is
known that the  moduli space of complex Enriques
surfaces is connected (at first, it was remarked at the end of
\cite{N0}; also see \cite{Ho}).
Thus, it is sufficient to prove this statement for one Enriques surface.

Let $\tau $ be an involution on the lattice $L$ which acts on $L$ like the
involution $\tau $ of Enriques surfaces. Then the lattice
$L_\tau =U\oplus U(2)\oplus E_8(2)$, where $E_8$ is an even unimodular
negative lattice of the rank 8. Let $c_1, c_2$ be a basis of $U$ such that
$(c_1)^2=(c_2)^2=0$ and $c_1\cdot c_2=1$. Let us consider an involution
$\sigma $ on $L$ such that $\sigma (c_1)=c_2, \sigma (c_2)=c_1$ and
$\sigma \mid (U)_L^\perp =-id$. The involution $\sigma $ evidently has
$L^\sigma =[c_1+c_2],$ where $(c_1+c_2)^2=2$. Thus, $L^\sigma $
is a hyperbolic lattice. We consider a  general element
$\omega _+\in L^\sigma \otimes \r$ and
$\omega _-\in L_{\tau, \sigma }\otimes \r $, where
$(\omega _+)^2=(\omega _-)^2>0$. By global
 Torelli theorem \cite{P\v S-\v S} and
epimorphicity of Torelli map \cite{Ku} for $K3$-surfaces, there exists a real
$K3$-surface $X$ with periods $\omega  _++i\omega _-$
such that $L=H^2(X(\cc );\z)$,
$\Pic X_\cc=L^\tau $,the action of the
antiholomorphic involution of $X$ is $\sigma $.  See \cite{N3} and
\cite{N4}
for details. Moreover, since $\Pic X_\cc=L^\tau $,
$X(\cc )$ has
a holomorphic involution which acts on $L$ as $\tau $. We denote this
automorphisms of $X$ by the same letters $\sigma $ and $\tau $
respectively. By
construction, the involution $\tau $ has no fixed points on $X(\cc )$
(see (3--3--4)) and by global Torelli Theorem
$\tau \sigma =\sigma \tau $
(since this is true for their action on cohomology).
 It follows that $X(\cc)/\{id, \tau \} $ is an Enriques surface with
antiholomorphic involution $\theta $ such that  $\sigma $ is
a lifting of $\theta $.

For $X(\cc )$, the group
${\Cal L}_{(\tau )}/{\Cal L}^{(\tau )}=
\{ c_1\mod 2, c_2\mod 2,c_1+c_2\mod 2\} .$ The subgroup
${\Cal F}/{\Cal L}^{(\tau )}$ should be $\sigma $-invariant, because
it is defined uniquely by the topology of $X(\cc )$ and $Y(\cc )$.

By construction of $\sigma $, we have
$\sigma (c_1 +L^\tau )=(c_2 +L^\tau )$,
$\sigma (c_2 +L^\tau )=(c_1 +L^\tau )$ and
$\sigma (c_1+c_2+L^\tau )=(c_1+c_2+L^\tau )$. Thus,
${\Cal F}/{\Cal L}^{(\tau )}=(c_1+c_2+L^\tau )$. We have
$(c_1+c_2)^2/2=1\mod 2$. Hence the Proposition follows.
\enddemo

\remark{Remark 3.4.1} It would be interesting to find a topological
proof of Lemma 3.4.3 if it does exist.
The Lemma is very important for complex Enriques surfaces.
For example, using this statement,
one can get more easily results of \cite{M-N}.
\endremark

We need the following statement too.

\proclaim{Lemma 3.4.4} Let $Y$ be a real Enriques surface and suppose
 there are two
liftings $\sigma $ and $\tau \sigma $ of the antiholomorphic involution
$\theta $ of $Y$ to antiholomorphic involutions of the surface
$X(\cc )$ (for example this is true if
$Y(\r )\not= \emptyset $). Then

$$
\dim H(\sigma )_+ +\dim H(\sigma )_-\le a(\sigma )\le
\dim H(\sigma )_+ +\dim H(\sigma )_- +1,
\tag3--4--7
$$
where the right inequality is a consequence of the following two
inequalities:
$$
r(\sigma )-a(\sigma )\ge -2\dim   H(\sigma )_++2r(\theta )
\tag3--4--8
$$
and
$$
r(\sigma )+a(\sigma )\le 2\dim  H(\sigma )_-+2r(\theta )+2.
\tag3--4--9
$$
\endproclaim

\demo{Proof} The inequality
$\dim H(\sigma )_+ +\dim H(\sigma )_-\le a(\sigma )$
follows from the general inequality of the condition 1.8.1 from
\cite{N4}. But we give a proof.

We use Notation 3.3.2.
By Lemma 3.3.1, $\Gamma _\pm =H(\sigma )_+\cap H(\sigma )_-$ is an isotropic
subgroup of $A_{q(\sigma )}$ and is orthogonal to
$H(\sigma )_+ + H(\sigma )_-$
with respect to the bilinear form of $q(\sigma )$. Thus,
with respect to the bilinear form of $q(\sigma )$ on $A_{q(\sigma )}$ we
have $H(\sigma )_+ + H(\sigma )_-
\subset \Gamma _\pm^{\ \perp} \subset A_{q(\sigma)}$,
where $\dim \Gamma_\pm^{\ \perp} =
\dim A_{q(\sigma )}-\dim \Gamma_\pm =a(\sigma )-
\dim \Gamma_\pm$,
since the bilinear form of $q(\sigma )$ is non-degenerate.
It follows that
$
\dim H(\sigma )_+ + \dim H(\sigma )_- - \dim \Gamma_\pm =
\dim (H(\sigma )_+ + H(\sigma )_-)\le \dim \Gamma_\pm^{\ \perp} =
a(\sigma )-\dim \Gamma_\pm$. Thus,
$\dim H(\sigma )_+ +\dim H(\sigma )_-\le a(\sigma )$.

The inequalities 1.8.2, 2) and 4) from the paper \cite{N4},
for our case give the
inequalities (3--4--8) and (3--4--9).
We mention  that these inequalities are equivalent to the
inequalities
$\dim  A_{L^\sigma _\tau }\le \rk  L^\sigma _\tau $ and
$\dim  A_{L_{\tau ,\sigma }}\le \rk L_{\tau ,\sigma }$
and one can prove them for this case independently.
Notation in \cite{N4} and here are connected as follows:
$t_{(+)}+t_{(-)}=r({\sigma }), a=a({\sigma }),
a_{H_{+}}=\dim H(\sigma )_+, a_{H_{-}}=\dim H(\sigma )_-,
p_{(+)}+p_{(-)}=r(\theta ), p_0=s_0=0,
l(A_{(S_+/C_+)^2})=r(\theta ), s_{(+)}+s_{(-)}=10,
l(A_{(S_-/C_-)^2})=10-r(\theta ), l_{(+)}+l_{(-)}=22.$

Multiplying the inequality (3--4--8) on $-1,$ we get the inequality
 $$
a(\sigma )-r(\sigma )\le 2\dim   H(\sigma )_+-2r(\theta ).
\tag3--4--10
$$
Adding the inequalities (3--4--9) and (3--4--10), we get the
right inequality (3--4--7).
\enddemo

Using Lemma 3.4.4, we can introduce a new invariant which is very important:

\definition{Definition 3.4.1} Let $Y$ be a real Enriques surface and
$\sigma $ is a lifting of the antiholomorphic involution $\theta$ of $Y$ to
an antiholomorphic involution of $X$. The invariant $\alpha (\sigma)$
is equal to
$$
\alpha (\sigma )=a(\sigma )-\dim H(\sigma )_+ -\dim H(\sigma )_- .
$$
By Lemma 3.4.4, $\alpha (\sigma )=0$ or $1$.
\enddefinition

By (3--4--1), the subspace
$$H^2(Y(\cc );\z )\otimes \z /2\subset H^2(Y(\cc );\z /2)$$
is a subspace of codimension 1. Thus, we can introduce the invariant:

\definition{Definition 3.4.2} Let $Y$ be a real Enriques surface with
an antiholomoprhic involution $\theta$. The invariant
$$
\beta (Y)=
\dim H^2(Y(\cc );\z /2)^\theta
-\dim (H^2(Y(\cc );\z )\otimes \z /2)^\theta .
$$
By (3--4--1), we have $0\le \beta (Y)\le 1$.
\enddefinition

Besides, we need some new invariants of
$\sigma $ from \cite{N4}.
In \cite{N4} for a triplet $(L,S,\sigma )$ (see \S 3.3) the invariants
$\delta _{\sigma S_+}$ and $\delta _{\sigma S_-}$ were defined. These
invariants are necessary for us for $S=L^\tau , S_+=L^{\tau ,\sigma },
S_-=L_\sigma ^\tau $.
We use Notation 3.3.2 to explain these invariants.

Let $A$ be a 2-elementary abelian group with a symmetric bilinear form
$b:A\times A\to (1/2)\z/\z$. An element $v\in A$ is called characteristic if
$b(x,x)=b(x,v)$ for any $x\in A$.
Since the
discriminant bilinear form $b(\sigma)$ (or bilinear form
of $q(\sigma)$) on the $2$-elementary
group $A_{q(\sigma)}$ is non-degenerate, there exists the unique
characteristic element
$v_{q(\sigma)}\in A_{q(\sigma)}$ for $b(\sigma)$.

The invariant $\delta _{\sigma L^{\tau ,\sigma }}=0$
if we simultaneously have:

(a) the characteristic element $v_{q(\sigma)}\in H(\sigma)_+$,

(b) The element $v_{q(\sigma )}\in H(\sigma )_+$ is
a characteristic element of the form
$b_{L^{\tau,\sigma}}\mid A_{L^{\tau,\sigma}}^{\ (1)}$
(see (3--3--13)).

Otherwise, the invariant $\delta _{\sigma L^{\tau ,\sigma }}=1$.

The definition of  $\delta_{\sigma L_\sigma ^\tau}$ is similar:
one should
consider $H(\sigma )_-$ and $L_\sigma ^\tau $ instead
$H(\sigma )_+$ and $L^{\tau,\sigma}$.

By Lemma 3.3.1, we have the following properties:
$$
\delta _{\sigma L^{\tau ,\sigma }}=0\ iff\
v_{q(\sigma)}\in \Gamma_\pm .
\tag3--4--11
$$
And
$$
\delta_{\sigma L_\sigma ^\tau}=0\ iff\
v_{q(\sigma)}\in \Gamma_\pm .
\tag3--4--12
$$
In particular, at any case,
$$
\delta _{\sigma L^{\tau ,\sigma }}=\delta_{\sigma L_\sigma ^\tau}.
\tag3--4--13
$$
By definition the invariant $\delta (\sigma )$ (see (3--3--3)),
the characteristic element
$v_{q(\sigma)}=0$ iff $\delta (\sigma )=0$. Thus, we have the
following important property too:
$$
\delta _{\sigma L^{\tau ,\sigma }}=\delta_{\sigma L_\sigma ^\tau}=0\ \
if\ \
\delta (\sigma )=0.
\tag3--4--14
$$

We mention another definition of these invariants. An element
$v(\sigma)\in L/2L$ is called characteristic if
$$
x \cdot \sigma (x)\equiv x\cdot v(\sigma )\mod 2.
$$
for any $x\in L$. The element $v(\sigma )$
is defined uniquely since $L$ is unimodular.
 From \cite{N4, \S 1} and the remark above, it follows that
$$
\delta _{\sigma L^{\tau ,\sigma }}=\delta_{\sigma L_\sigma ^\tau}=0\ \
iff\ \ v(\sigma )\in
L^{\tau ,\sigma } \cap L_\sigma ^\tau \mod 2L
\tag3--4--15
$$
The geometrical sense of $v(\sigma )$ is that
$$
X_\sigma(\r ) \sim v(\sigma ) \mod 2\ in\ H_2(X(\cc );\z).
\tag3--4--16
$$

Now we can formulate the basic Theorem.

\proclaim
{Theorem 3.4.5} Let $Y$ be a real Enriques surface and suppose
 there are two
liftings $\sigma $ and $\tau \sigma $ of the antiholomorphic involution
$\theta $ of $Y$ to antiholomorphic involutions of the $K3$-surface
$X(\cc )$ (for example this is true if
$Y(\r )\not= \emptyset $).

Then
$$
\dim H^2(Y(\cc );\z /2)^\theta =
\dim (H^2(Y(\cc );\z )\otimes \z /2)^\theta =10-a(\theta )\ and\
\beta(Y)=0
$$
if
$\alpha(\sigma )=1$ and
$\delta _{\sigma L^{\tau ,\sigma }}=\delta _{\sigma L^\tau _\sigma }=0$.

Otherwise (if
either $\alpha(\sigma )=0$
or $\delta _{\sigma L^{\tau ,\sigma }}+\delta _{\sigma L^\tau _\sigma }>0$),
$$
\dim (H^2(Y(\cc );\z )\otimes \z /2)^\theta =11-a(\theta ),
$$
and
$$
\dim H^2(Y(\cc );\z /2)^\theta =11-a(\theta )+\beta (Y),\ \
where \ \ 0\le \beta (Y)\le 1.
$$
\endproclaim

\demo{Proof}  All homomorphisms (3--4--1) ---(3--4--6) are equivariant
with respect to the action of $\theta $ and its lifting $\sigma $.
Thus, by (3--4--6)
$$
(H^2(Y(\cc );\z )\otimes \z /2)^\theta \cong ({\Cal L}/{\Cal F})^\sigma
$$
and
we should calculate
$$\dim ({\Cal L}/{\Cal F})^\sigma .$$

The lattice $L$ is  unimodular  and defines
the natural exact sequence
$$
0\to L^\tau \to L \to L^\ast _\tau \to 0
$$
 and the corresponding exact sequence
$$
0\to {\Cal L}^{(\tau )} \to {\Cal L} \to {\Cal L}^\ast _{(\tau )}\to 0.
$$

Thus, we have the canonical isomorphism
${\Cal L}/{\Cal L}^{(\tau )}\cong {\Cal L}^\ast _{(\tau )}$ and
the corresponding filtration of subgroups
$$
F={\Cal F}/{\Cal L}^{(\tau )}\subset H={\Cal L}_{(\tau )}/{\Cal L}^{(\tau )}
\subset {\Cal L}^\ast _{(\tau )},
$$
where $\dim F=1$ and $\dim H=2$.
We should calculate $\dim ({\Cal L}^\ast _{(\tau )}/F)^\sigma $. Let
$F=\{0 ,f\} $, where $f$ is some
linear function on ${\Cal L}_{(\tau )}$. Let
$\phi \in {\Cal L}^\ast _{(\tau )}$. Then $\phi +F\in
({\Cal L}^\ast _{(\tau )}/F)^\sigma $ iff $\sigma ^\ast \phi +\phi \in F$.
Thus, for $\phi $ we have two possibilities.

Case 1. $\sigma ^\ast \phi +\phi =0$
(or $\phi \in ({\Cal L}^\ast _{(\tau )})^\sigma $). Thus, equivalently
$\phi ((id+\sigma )({\Cal L}_{(\tau )}))=0$.
 From the exact sequence
$$
0\to ({\Cal L}_{(\tau )})^\sigma
\to {\Cal L}_{(\tau )}@>id+\sigma >>{\Cal L}_{(\tau )}
$$
we then have that $\dim ({\Cal L}^\ast _{(\tau )})^\sigma =
\dim ({\Cal L}_{(\tau )})^\sigma $.

Case 2. $\sigma ^\ast \phi +\phi =f$. Obviously, such a $\phi $ does
exist  iff $f\in Im(id +\sigma )^\ast $. Equivalently,  by the exact
sequence above,
$f(({\Cal L}_{(\tau )})^\sigma )=0.$

Thus, by Proposition 3.4.2,
 $$
\dim ({\Cal L}/{\Cal F})^\sigma = 11-a(\theta )-a(\sigma )+
\dim H(\sigma )_++\dim H(\sigma )_-
$$
if  $f(({\Cal L}_{(\tau )})^\sigma )\not=0$, and
$$
\dim ({\Cal L}/{\Cal F})^\sigma = 12-a(\theta )-a(\sigma )+
\dim H(\sigma )_++\dim H(\sigma )_-
$$
 if  $f(({\Cal L}_{(\tau )})^\sigma )=0$. Then the statements of
the Theorem about
$\dim (H^2(Y(\cc );\z )\otimes \z /2)^\theta $ follow from

\proclaim{Lemma 3.4.6}
 If $\alpha(\sigma )=0$, then $f(({\Cal L}_{(\tau )})^\sigma )\not=0$.

 If
$\alpha(\sigma ) = 1$ and
$\delta _{\sigma L^{\tau ,\sigma }}+\delta _{\sigma L^\tau _\sigma }>0$,
then $f(({\Cal L}_{(\tau )})^\sigma )=0$.

If
$\alpha (\sigma ) = 1$ and
$\delta _{\sigma L^{\tau ,\sigma }}=\delta _{\sigma L^\tau _\sigma }=0$,
then $f(({\Cal L}_{(\tau )})^\sigma )\not=0$.
\endproclaim

\demo{Proof}
 From
\cite{N4, \S 1, formula (8.9)}, it follows that
$$
\rk L^\sigma _\tau -\rk A_{L^\sigma _\tau }=(r(\sigma )-a(\sigma ))-
(-2\dim H(\sigma )_++2r(\theta ))
\tag3--4--17
$$
and
$$
\rk  L_{\tau ,\sigma }-\rk A_{L_{\tau ,\sigma }}=
(2\dim H(\sigma )_-+2r(\theta )+2)-(r(\sigma )+a(\sigma )).
\tag3--4--18
$$

Let $a(\sigma )=\dim H(\sigma )_++\dim H(\sigma )_-+1$ (equivalently,
$\alpha (\sigma )=1$). Then
both inequalities (3--4--8) and (3--4--9)  are equalities
(see the proof of Lemma 3.4.4)
and from
(3--4--17) and (3--4--18) we then get that
$\rk L^\sigma _\tau =\rk A_{L^\sigma _\tau }$ and
$\rk  L_{\tau ,\sigma }=\rk A_{L_{\tau ,\sigma }}$.

By Proposition 3.4.2,
$\dim ({\Cal L}^{(\tau )})^\sigma =10-a(\theta ),
\dim ({\Cal L}_{(\tau )})^\sigma =11-a(\theta ).$ It follows that
$$
({\Cal L}_{(\tau )})^\sigma +{\Cal L}^{(\tau )}\subset
{\Cal L}_{(\tau )}/{\Cal L}^{(\tau )}=H
$$
is a subspace of dimension 1 (over $\z /2$). By Proposition 3.4.2,
$({\Cal L}_{(\tau )})^\sigma =
(L^\sigma _\tau \oplus L_{\tau ,\sigma })/2L_\tau $. From Lemma 3.4.3
and earlier considerations, we then get that \linebreak
$f(({\Cal L}_{(\tau )})^\sigma )\not=  0$ iff
$$
(L_\tau ^\sigma )^2/2\subset 2\z  \ and \ \
(L_{\tau ,\sigma })^2/2\subset 2\z .
\tag 3--4--19
$$
Since
$\rk L^\sigma _\tau =\rk A_{L^\sigma _\tau }$ and
$\rk  L_{\tau ,\sigma }=\rk A_{L_{\tau ,\sigma }}$ and from the
decomposition of 2-adic lattices as an orthogonal sum of
elementary lattices of rank 1 or 2, we get that (3--4--19) is
equivalent to the facts
$$
q_{L^\sigma _\tau }\ne  q_{\gamma _1}^{(2)}(2)\oplus q_1
\ \ \ (\gamma _1\in \z _2^\ast /(\z _2^\ast )^2),
\tag3--4--20
$$
and
$$
q_{L_{\tau ,\sigma }}\ne q_{\gamma _2}^{(2)}(2)\oplus q_2
\ \ \ (\gamma _2\in \z _2^\ast /(\z _2^\ast )^2),
\tag3--4--21
$$
where $q^{(2)}_\gamma (2)$
is a non-degenerate quadratic form on a group of order 2. From
\cite{N4, \S 1, (8.8)---(8.10)}, it follows that (3--4--20) is equivalent to
$\delta _{\sigma L^{\tau, \sigma }}=0$ and (3--4--21) is equivalent to
$\delta _{\sigma L^\tau _\sigma }=0$. Thus,
we get the statement of Lemma 3.4.6 for
$\alpha (\sigma )=1$.

Now let $a(\sigma )=\dim H(\sigma )_++\dim H(\sigma )_-$
(equivalently, $\alpha (\sigma )=0$). Then, from
Lemma 3.4.4,  (3--4--17) and (3--4--18), we get that either
$2\ge \rk L^\sigma _\tau -\rk A_{L^\sigma _\tau }\ge 1$ or
$2\ge \rk  L_{\tau ,\sigma }-\rk A_{L_{\tau ,\sigma }}\ge 1$.
For example, let $\rk L^\sigma _\tau -\rk A_{L^\sigma _\tau }\ge 1.$
Since the lattice $L^\sigma _\tau$  is even,
$\rk L^\sigma _\tau =\rk A_{L^\sigma _\tau }+2.$ It follows that
$L^\sigma _\tau \otimes \z _2\cong K\oplus M$, where $K$ is an
unimodular 2-adic lattice of the rank 2. Then
$L_\tau ^\sigma \mod 2L_\tau  ={\Cal L}_{(\tau )}/{\Cal L}^{(\tau )}=H,$
since ${\Cal L}^{(\tau )}$ is a kernel of the form on ${\Cal L}_{(\tau )}$
and $H$ has the dimension 2. Since
$({\Cal L}_{(\tau )})^\sigma =
(L_\tau ^\sigma \oplus L_{\tau ,\sigma })/2L_\tau  $, we get
 $f(({\Cal L}_{(\tau )})^\sigma )\not= 0$, since
$({\Cal L}_{(\tau )})^\sigma \mod {\Cal L}^{(\tau )}$ contains elements
$g_1,g_2$ (see considerations before Lemma 3.4.3) and
$f\cdot g_1=f\cdot g_2=1\mod 2$.

To finish the proof of Theorem 3.4.5, we should now show that
$\beta (Y)=0$ if simultaneously
$a(\sigma )=\dim H(\sigma )_++\dim H(\sigma )_-+1$
(or $\alpha (\sigma )=1$) and
$\delta _{\sigma L^{\tau, \sigma }}=\delta _{\sigma L^\tau _\sigma }=0$.

 From (3--4--4),
$(H^2(Y(\cc );\z )\otimes \z /2)/\Tor H^2(Y(\cc );\z )\cong
{\Cal L}^{(\tau )}$. From Proposition 3.4.2, we then get that
$$
\dim ((H^2(Y(\cc );\z )\otimes \z /2)/\Tor H^2(Y(\cc );\z ))^\theta =
10-a(\theta ).
$$
Besides, we had proved that
$$
\dim (H^2(Y(\cc );\z )\otimes \z /2)^\theta
=10-a(\theta ).
$$
It follows that
$$
\dim H^2(Y(\cc );\z )\otimes \z /2)^\theta
=\dim ((H^2(Y(\cc );\z )\otimes \z /2)/\Tor H^2(Y(\cc );\z ))^\theta .
$$
We now show that from this equality, we have
$$
H^2(Y(\cc );\z /2)^\theta =
(H^2(Y(\cc );\z )\otimes \z /2)^\theta .
$$
Let $0\not=t\in \Tor H^2(Y(\cc );\z ).$ Since $t$ is a
torsion element,
the intersection pairing $t \cdot x=0 $ for any integral class
$x\in H^2(Y(\cc );\z )$,
because $t\cdot x\in H^4(Y(\cc );\z )\cong \z $. Thus,
$t\perp  H^2(Y(\cc );\z )\otimes \z /2$ with respect to the
intersection
pairing. Since $H^2(Y(\cc );\z )\otimes \z /2\subset
H^2(Y(\cc );\z /2)$ has codimension 1 and the intersection pairing on
$H^2(Y(\cc );\z /2)$ is non-degenerate, it follows that
$t\cdot u=1\mod 2$ for any
$$
u\in H^2(Y(\cc );\z /2)-H^2(Y(\cc );\z )\otimes \z /2.
$$
Suppose that
$$
\beta (Y)=\dim H^2(Y(\cc );\z /2)^\theta
-\dim (H^2(Y(\cc );\z )\otimes \z /2)^\theta >0.
$$
Then, there exists
$$
u\in H^2(Y(\cc );\z /2)^\theta -(H^2(Y(\cc );\z )\otimes \z /2)^\theta .
$$
We have seen that $u\cdot t\not=0\mod 2$. Since $\theta (u)=u$
and $\dim  \Tor H^2(Y(\cc );\z )=1$,
we then get the $\theta $-invariant decomposition
$$
H^2(Y(\cc );\z )\otimes \z /2=\Tor H^2(Y(\cc );\z )\oplus
\Delta
$$
where
$t\cdot \Delta =0\mod 2 $. Also, we have
$$
\dim (H^2(Y(\cc );\z )\otimes \z /2)^\theta =1+\dim \Delta ^\theta ,
$$
where $\Delta \cong (H^2(Y(\cc );\z )\otimes \z /2)/\Tor H^2(Y(\cc );\z )$.
Thus,
$$
\dim \Delta ^\theta =
\dim ((H^2(Y(\cc );\z )\otimes \z /2)/\Tor H^2(Y(\cc );\z ))^\theta .
$$
Thus, we get that
$$
\dim (H^2(Y(\cc );\z )\otimes \z /2)^\theta
=1+\dim ((H^2(Y(\cc );\z )\otimes \z /2)/\Tor H^2(Y(\cc );\z ))^\theta .
$$
This gives rise to a contradiction and
finishes the proof of Theorem 3.4.5.
\enddemo
\enddemo

 From Proposition 3.4.1 and Theorem 3.4.5,
we get the following basic result of this section (see
Definitions 1.2 and 3.4.2):

\proclaim{Theorem 3.4.7}
Let $Y$ be a real Enriques surface such that there are two
liftings $\sigma $ and $\tau \sigma $ of the antiholomorphic involution
$\theta $ of $Y$ to antiholomorphic involutions of the $K3$-surface
$X(\cc )$ (this is true for example if
$Y(\r )\not= \emptyset $).

Then
$$
b(Y)=r(\theta )-a(\theta )\ge 0\ and\ \beta(Y)=0
$$
if $\alpha(\sigma )=1$ and
$\delta _{\sigma L^{\tau ,\sigma }}=\delta _{\sigma L^\tau _\sigma }=0$.

Otherwise (if we have
either $\alpha(\sigma )=0$ or
$\delta _{\sigma L^{\tau ,\sigma }}+\delta _{\sigma L^\tau _\sigma }>0$)
$$
b(Y)=r(\theta )-a(\theta )+1+\beta (Y)\ge 1+\beta (Y).
$$
where $\beta(Y)=0$ or $1$.
\endproclaim

We will remark that inequalities here follow from the
evident inequality $r(\theta )\ge a(\theta )$. This holds since
$$
r(\theta )=\rk L^{\tau ,\sigma }
\ge \rk (L^{\tau ,\sigma }(1/2))^\ast /L^{\tau ,\sigma }(1/2)=a(\theta ).
$$

We can unify the formulae of Theorem 3.4.7:
$$
b(Y)=r(\theta)-a(\theta)+
\max \{ 1 - \alpha(\sigma),
(\delta _{\sigma L^{\tau ,\sigma }} + \delta _{\sigma L^\tau _\sigma })/2 \}+
\beta(Y)
\tag3--4--22
$$
where $\beta(Y)=0$ if
$\alpha(\sigma)=1$ and
$\delta _{\sigma L^{\tau ,\sigma }} = \delta _{\sigma L^\tau _\sigma }=0$.
We recall that by (3--4--13),
$\delta _{\sigma L^{\tau ,\sigma }} = \delta _{\sigma L^\tau _\sigma }$.

\subhead
3.5. Connection between $b(Y)$ and
$s(\sigma )+s(\tau \sigma )=s_{nor}+2s_{or}$
\endsubhead

 From Theorems 3.3.4, 3.4.7 and the equality
$a(\sigma )=\dim H(\sigma )_+ + \dim H(\sigma )_- +
\alpha (\sigma )$, we get the following basic statements of
this paper. We divide them into three parts according to
the cases of Theorem 3.4.7.

\proclaim{Theorem 3.5.1}
Let $Y$ be a real Enriques surface such that there are two
liftings $\sigma $ and $\tau \sigma $ of the antiholomorphic involution
$\theta $ of $Y$ to antiholomorphic involutions of the $K3$-surface
$X(\cc )$ (this is true for example if
$Y(\r )\not= \emptyset $). Suppose that  $\alpha (\sigma )=1$ and
$\delta _{\sigma L^{\tau ,\sigma }}=\delta _{\sigma L^\tau _\sigma }=0$.

Then $\beta (Y)=0$ and:

If both $s(\sigma )>0$ and $s(\tau \sigma )>0$,
$$
b(Y)=s(\sigma )+s(\tau \sigma )-2+\dim H(\sigma )_- -
\dim (H(\sigma )_+)^\perp \cap H(\sigma )_-.
$$

If either $s(\sigma )=0, s(\tau \sigma )>0$ or
$s(\tau \sigma )=0, s(\sigma )>0$,
$$
b(Y) = s(\sigma ) + s(\tau \sigma ) - 1 + \dim H(\sigma )_- -
\dim (H(\sigma )_+)^\perp \cap H(\sigma )_-.
$$

If $s(\sigma )=s(\tau \sigma )=0$,
$$
b(Y) = s(\sigma )+s(\tau \sigma )+\dim H(\sigma )_- -
\dim  (H(\sigma )_+)^\perp
\cap H(\sigma )_-.
$$

Here $s(\sigma )+s(\tau \sigma )=s_{nor}+2s_{or}$.
\endproclaim

\proclaim{Theorem 3.5.2}
Let $Y$ be a real Enriques surface such that there are two
liftings $\sigma $ and $\tau \sigma $ of the antiholomorphic involution
$\theta $ of $Y$ to antiholomorphic involutions of the $K3$-surface
$X(\cc )$ (this is true for example if
$Y(\r )\not= \emptyset $). Suppose that  simultaneously
$\alpha (\sigma )=1$ and
$\delta _{\sigma L^{\tau ,\sigma }}+\delta _{\sigma L^\tau _\sigma }>0$.

Then both $s(\sigma )>0, s(\tau \sigma )>0$, and
$$
b(Y) = s(\sigma ) + s(\tau \sigma )-1 + \dim H(\sigma )_- -
\dim (H(\sigma )_+ )^\perp \cap H(\sigma )_- + \beta (Y).
$$

Here $s(\sigma )+s(\tau \sigma )=s_{nor}+2s_{or}$ and
$0\le \beta (Y)\le 1$.
\endproclaim

\proclaim{Theorem 3.5.3}
Let $Y$ be a real Enriques surface such that there are two
liftings $\sigma $ and $\tau \sigma $ of the antiholomorphic involution
$\theta $ of $Y$ to antiholomorphic involutions of the $K3$-surface
$X(\cc )$ (this is true for example if
$Y(\r )\not= \emptyset $). Suppose that
$\alpha(\sigma )=0$.

Then:

If both $s(\sigma )>0$ and $s(\tau \sigma )>0$,
$$
b(Y)=s(\sigma )+s(\tau \sigma )-2 + \dim H(\sigma )_- -
\dim  (H(\sigma )_+)^\perp \cap H(\sigma )_-+\beta (Y)
\ge 1+\beta (Y).
$$

If either $s(\sigma )=0, s(\tau \sigma )>0$ or
$s(\tau \sigma )=0, s(\sigma )>0$,
$$
b(Y)=s(\sigma )+s(\tau \sigma )-1 + \dim H(\sigma )_- -
\dim  (H(\sigma )_+)^\perp
\cap H(\sigma )_- + \beta (Y) \ge 1+\beta (Y).
$$

If $s(\sigma )=s(\tau \sigma )=0,$
$$
b(Y) = s(\sigma )+s(\tau \sigma ) + \dim H(\sigma )_- -
\dim (H(\sigma )_+)^\perp
\cap H(\sigma )_- + \beta (Y) \ge 1+\beta (Y).
$$

Here $s(\sigma )+s(\tau \sigma )=s_{nor}+2s_{or}$ and
$0\le \beta (Y)\le 1$.
\endproclaim

Similar to (3--4--22), we can write down an unifying formula:
$$
 b(Y) = s(\sigma)+s(\tau\sigma)-
\sharp \{x\in \{s(\sigma),s(\tau\sigma)\} \mid x>0 \}+
\tag3--5--1
$$
$$
+ \min \{ \alpha(\sigma),
(\delta _{\sigma L^{\tau ,\sigma }} + \delta _{\sigma L^\tau _\sigma })/2 \}
+ \dim H(\sigma )_- - \dim (H(\sigma )_+)^\perp
\cap H(\sigma )_- + \beta (Y)
$$
where $\beta(Y)=0$ if
$\alpha(\sigma)=1$ and
$\delta _{\sigma L^{\tau ,\sigma }} = \delta _{\sigma L^\tau _\sigma }=0$.

By these Theorems, it is important to estimate the invariant
$$
\gamma (\sigma )=\dim H(\sigma )_- -
\dim (H(\sigma )_+)^\perp \cap H(\sigma )_- .
\tag3--5--2
$$
Here we have
\proclaim{Proposition 3.5.4} We have inequalities:
$0\le \gamma (\sigma )\le 2$.
\endproclaim

\demo{Proof} We use Notation 3.3.2. By (3--3--10), either
$\dim H(\sigma )_+ - \dim \Gamma_\pm \le 2$ or
$\dim H(\sigma )_- - \dim \Gamma_\pm \le 2$
because
$\dim H(\sigma )_+\le r(\theta )$ and
$\dim H(\sigma )_-\le 10-r(\theta)$.
By Lemma 3.3.1,
$\Gamma _\pm $ is an isotropic subgroup of $A_{q(\sigma )}$ for
a bilinear form of $q(\sigma )$. It follows the Proposition.
\enddemo

We mention another inequalities for items of the
formula (3--5--1): Evidently,
$$
2\ge \sharp \{x\in \{s(\sigma),s(\tau\sigma)\} \mid x>0 \} \ge 0,
\tag3--5--3
$$
$$
1 \ge \min \{ \alpha(\sigma),
(\delta _{\sigma L^{\tau ,\sigma }} + \delta _{\sigma L^\tau _\sigma })/2 \}
\ge 0.
\tag3--5--4
$$
Besides, we have
$$
\sharp \{x\in \{s(\sigma),s(\tau\sigma)\} \mid x>0 \} -
\min \{ \alpha(\sigma),
(\delta _{\sigma L^{\tau ,\sigma }} + \delta _{\sigma L^\tau _\sigma })/2 \}
\ge 0,
\tag3--5--5
$$
and we have an equality here only if $s(\sigma )=s(\tau\sigma )=0$.
We explain the last inequality.
If $s(\sigma)=0$, then $\delta(\sigma )=0$
and
$\delta _{\sigma L^{\tau ,\sigma }} = \delta _{\sigma L^\tau _\sigma } = 0$
(see (3--3--5) and (3--4--14)). It follows the (3--5--5).

Further, we give applications of our basic Theorems
1.2, 2.1, 3.4.7 and 3.5.1---3.5.3.

\subhead
3.6. The case $Y(\r )\not=\emptyset$ and $b(Y)=0$
\endsubhead

This case is interesting because, by Theorem 1.2 and Proposition 0.3,  in
this case
we know precisely the dimension:
$\dim~_2Br(Y)=\epsilon (Y)=1$.

The following
Theorem describes completely this situation.

\proclaim{Theorem 3.6.1} Let $Y$ be a real Enriques surface such that
$Y(\r )\not=\emptyset$.

Then $b(Y)=0$ iff $Y(\r )$
is connected non-orientable and the invariants $r(\theta )$ and
$a(\theta)$ are equal: $r(\theta )=a(\theta )$.

For this surface $Y$, the invariant $\epsilon (Y)=1$ and
$\dim \ _2Br(Y)=1$.
\endproclaim

\demo{Proof} We use Notation 3.3.2.
Let $b(Y)=0$. Since $Y(\r )\not=\emptyset$,
by Theorem 2.1, we have $s=1$. Thus $Y(\r )$ is
connected. By Lemma 3.2.1,
one of the involutions $\sigma $ or $\tau \sigma $ has empty set of real
points. We can suppose that $s(\sigma )=0$. By (3--3--5) and
(3--4--14), then the invariants $\delta (\sigma )=0$ and
$\delta _{\sigma L^{\tau ,\sigma }}=\delta _{\sigma L^\tau _\sigma }=0$.
By Theorems  3.5.1, 3.5.2 and 3.5.3, we then get that
$\alpha (\sigma )=1$, $s(\tau \sigma )=1$ and
$H(\sigma )_+ \perp H(\sigma )_-$.
Since $s(\sigma )=0$ and $s(\tau \sigma )=1$, by Lemma 3.2.1,
the surface $Y(\r )$ is non-orientable.
Since $b(Y)=0$, from Theorem 3.4.7,
$r(\theta )=a(\theta )$.

Now let $Y(\r )$ be a connected non-orientable surface and
$r(\theta )=a(\theta )$.
By Lemma 3.2.1,
$$
s(\sigma )+s(\tau \sigma )=1.
\tag3--6--1
$$
We can suppose that
$s(\sigma )=0$ and $s(\tau \sigma )=1$.

We claim that  if $r(\theta )=a(\theta )$, then
$H(\sigma )_+ \perp H(\sigma )_-$.
Actually, we have a subgroup $\Gamma (\sigma )_\pm \subset H(\sigma )_+$
where
$\dim \Gamma (\sigma )_\pm = a(\theta )$ and
$\dim H(\sigma )_+\le r(\theta )$. Since $r(\theta )=a(\theta )$, we
then get that $\Gamma (\sigma )_\pm = H(\sigma )_+$.
Thus, we have
$\Gamma (\sigma )_\pm = H(\sigma )_+ \subset H(\sigma )_-$.
By Lemma 3.3.1,
the subgroup $\Gamma (\sigma )_\pm$ is isotropic for the bilinear form
$q(\sigma )$ on $A_{q(\sigma )}$. Hence the
statement we claimed follows. Then we get that
$$
\gamma (\sigma )=\dim H(\sigma )_- -
\dim (H(\sigma )_+)^\perp
\cap H(\sigma )_-=0.
\tag3--6--2
$$
Like above,  $\delta (\sigma )=0$ and
$\delta _{\sigma L^{\tau ,\sigma }}=\delta _{\sigma L^\tau _\sigma }=0$.
Thus, we can apply to our situation Theorems 3.5.1 or 3.5.3. If
$\alpha(\sigma )=0$, we should apply
Theorem 3.5.3. Then from (3--6--1) and (3--6--2) we get the
contradiction:
$0=b(Y)\ge 1+\beta (Y)>0$.
Thus, $\alpha (\sigma)=1$, and we should
apply Theorem 3.5.1. Then, from (3--6--1) and (3--6--2), we get
$b(Y)=0$.

By Theorem 1.2 and Proposition 0.3, we have
$\dim \ _2Br(Y)=\epsilon (Y)=1$, since $b(Y)=0$ and
$Y(\r )\not= \emptyset $.
This proves the Theorem.
\enddemo

\subhead
3.7. Some geometric applications
\endsubhead

We want to prove the following result
where we use the invariants described above.

\proclaim{Theorem 3.7.1} Let $Y$ be a real Enriques surface.
Then we have the inequalities for the numbers $s$ and $s_{nor}$
or connected and non-orientable connected components of $Y(\r)$:
$$
s\le (2+r(\theta )-a(\theta)+\max \{ 1 - \alpha(\sigma),
(\delta _{\sigma L^{\tau ,\sigma }} + \delta _{\sigma L^\tau _\sigma })/2 \}+
\beta (Y))/2 \le 6;
$$
$$
\split
s_{nor}
&\le 2-  \sharp \{x\in \{s(\sigma),s(\tau\sigma)\} \mid x>0 \}+
\min \{ \alpha(\sigma),
(\delta _{\sigma L^{\tau ,\sigma }} + \delta _{\sigma L^\tau _\sigma })/2 \}\\
&+ \dim H(\sigma )_- - \dim (H(\sigma )_+)^\perp
\cap H(\sigma )_- + \beta (Y) \le 4.
\endsplit
$$
\endproclaim

\demo{Proof} The first inequality follows from
(3--4--22), Theorem 2.1 and (3--3--10). The second inequality follows from
the formula (3--5--1), Theorem 2.1, Proposition 3.5.4 and
inequalities (3--5--2)---(3--5--5).
\enddemo

We should mention that the inequality
$s\le 6$ also follows from results of V. M. Harlamov
\cite{Ha2}

\subhead
3.8. A remark about further results
\endsubhead

We want to mention further results which were obtained by
the first author during the time this paper was considered for a
publication.

In \cite{N5}, further results about Brauer groups of
real Enriques surfaces were obtained. Many of them valid for an
arbitrary real smooth projective
algebraic surface. We cite these results below.

The following Theorem generalizes the result from \cite{C-P} we
mentioned in Introduction about the homomorphism
(0--1).

\proclaim{Theorem 3.8.1} Let $X/\r $ be an algebraic projective manifold
with the antiholomorphic involution $g$, and $G=\{ \id, g\} $.
Then the homomorphism (0--1) is epimorphic if
$H^3(X(\cc )/G;\z/2)=0$. More generally, the homomorphism (0--1) is
epimorphic if the kernel of the homomorphism
$$
i^*:H^3(X(\cc )/G;\z/2)\to H^3(X(\r );\z/2)
$$
is equal to zero. Here $i:X(\r )\subset X(\cc )/G$ denote the
embedding.
\endproclaim

The following result valid for smooth surfaces.

\proclaim{Theorem 3.8.2}
Let $X$ be  a real smooth projective algebraic surface and \linebreak
$H^3(X(\cc )/G;\z/2)=0$.

Then
the Hochschild--Serre spectral sequence
degenerates and
$$
\dim~_2Br(X)=2s-1+h^{2,0}(X(\cc ))+h^{1,1}_-(X(\cc ))-
\rho_ +(X\otimes \cc ).
$$
Here $h_-^{1,1}(X(\cc ))=\dim H_-^{1,1}(X(\cc ))$ where
$$
H_-^{1,1}(X(\cc ))=\{x\in H^{1,1}(X(\cc ))\mid g(x)=-x \}
$$
is the set of potentially real algebraic cycles of $X$. And
$\rho _+(X\otimes \cc )=\dim (\Pic (X\otimes \cc)\otimes \cc)^G$.
The characteristic cycle map gives an injection of
$(\Pic (X\otimes \cc)\otimes \cc)^G$ to $H_-^{1,1}(X(\cc ))$.
\endproclaim

Both these Theorems work for "general" real
Enriques surfaces, and we get for these surfaces exactly the same
results as for rational surfaces with non-empty set of real points.

\proclaim{Theorem 3.8.3} Let $Y$ be a real Enriques surface.

Then  $H^3(Y(\cc )/G;\z/2)=0$ iff
the both liftings $\sigma $ and $\tau\sigma $ of
the antiholomorphic involution $\theta $ of $Y$ have non-empty
sets $X_\sigma (\r )$ and $X_{\tau\sigma }(\r )$ of real points
($s(\sigma )>0$ and $s(\tau\sigma )>0$).

Thus (by Theorems 3.8.1, 3.8.2),
in this case, the homomorphism (0--1) is epimorphic,
Hochschild--Serre spectral sequence degenerates, and
$$
b(Y)=2s-2,\ \ \dim~_2Br(Y)=2s-1.
$$
\endproclaim

Now we want to show how Theorems 3.4.7 and 3.5.1---3.5.3 of our
paper work together with
Theorem 3.8.3.
By Theorem 3.4.7, if
$\alpha(\sigma )=1$ and
$\delta _{\sigma L^{\tau, \sigma}}=\delta _{\sigma L^\tau _\sigma}=0$,
then at any case, the invariant $\beta (Y)=0$.

Now suppose that either
$\alpha(\sigma )=0$ or
$\delta _{\sigma L^{\tau, \sigma}} + \delta _{\sigma L^\tau _\sigma}>0$.
Then, by Theorem 3.4.7,
$b(Y)=r(\theta )-a(\theta )+ 1 +\beta (Y)$. Additionally, let us suppose that
$s(\sigma )>0$ and $s(\tau\sigma )>0$. By Theorem 3.8.3, we
then get that
$b(Y)=2s-2\equiv 0\mod 2$. Since $r(\theta)-a(\theta )\equiv 0\mod 2$,
we have $\beta(Y)=1$. Thus, we have

\proclaim{Theorem 3.8.4} Let $Y$ be a real Enriques surface and
either
$\alpha(\sigma )=0$ or
$\delta _{\sigma L^{\tau, \sigma}}+\delta _{\sigma L^\tau _\sigma}>0$.

Then if $s(\sigma )>0$ and $s(\tau\sigma )>0$, we have
$\beta (Y)=1$.
\endproclaim

We don't know now, does the same statement holds if one of
$s(\sigma )$ or $s(\tau\sigma )$ is 0.

We continue to suppose that
$s(\sigma )>0, s(\tau\sigma )>0$. For example, assume the case
$\alpha(\sigma )=0$. Then by
Theorem 3.5.3 and Theorems 3.8.3 and 3.8.4, we have the equality:
$s_{nor}+2s_{or}-2+\dim~H(\sigma )_- -
\dim~(H(\sigma )_+)^\perp \cap H(\sigma )_- + 1=
2s_{nor}+2s_{or}-2$.
Thus, we have
$$
\split
& if\   \alpha (\sigma )=0,\ then \\
& s_{nor}=1+\dim H(\sigma )_- -
\dim~(H(\sigma )_+)^\perp \cap H(\sigma )_-.
\endsplit
\tag3--8--1
$$
By the same considerations, from
Theorems 3.8.3, 3.8.4 and Theorems 3.5.1, 3.5.2, we get
$$
\split
& if\
(\alpha (\sigma )=1,\
\delta _{\sigma L^{\tau, \sigma}}=\delta _{\sigma L^\tau _\sigma}=0),\
then \\
& s_{nor}=\dim~H(\sigma )_- -
\dim~(H(\sigma )_+)^\perp \cap H(\sigma )_-;
\endsplit
\tag3--8--2
$$
$$
\split
& if\
(\alpha (\sigma )=1,\
\delta _{\sigma L^{\tau, \sigma}}+\delta _{\sigma L^\tau _\sigma}>0),\ then \\
& s_{nor}=2 + \dim~H(\sigma )_- -
\dim~(H(\sigma )_+)^\perp \cap H(\sigma )_-;
\endsplit
\tag3--8--3
$$
Statements (3--8--1), (3--8--2) and (3--8--3)
are very important for a topological classification of
real Enriques surfaces (see below). We recall that we claim these
statements only if $s(\sigma )>0$ and
$s(\tau\sigma )>0$.

We remark that Enriques surfaces of Theorem 3.6.1 does not satisfy to
conditions of Theorem 3.8.3. For these surfaces one of
$s(\sigma ),\ s(\tau\sigma )$ is zero.
But there exists another generalization of the Theorem 3.6.1.

\proclaim{Theorem 3.8.5} Let $Y$ be a real Enriques surface.

Then the inequality (of Theorem 2.1)
$$
b(Y)\ge 2s-2
$$
is an equality iff the Hochschild--Serre spectral sequence degenerates.
In particular (by Theorem 1.2), we have
$$
\epsilon (Y)=1\ and\ \dim~_2Br(Y)=2s-1,\ if\ b(Y)=2s-2.
$$
At any case, we have an inequality
$$
\dim~_2Br(Y)\ge 2s-1.
$$

\endproclaim

We mention that the inequality $\dim~_2Br(Y)\ge 2s-1$
also gives another proof of the Corollary 2.2.

\smallpagebreak

In \cite{N6}, there were obtained some results about a topological
classification of real Enriques surfaces. We cite some of these results
below.

All possible invariants of \cite{N4} for triplets
$(L,\sigma, S)$ (see \S 3.3)
were described. In particular, all possibilities for pairs
$(X_\sigma (\r ), X_{\tau\sigma }(\r ))$ were obtained. Using formulae
(3--8--1), (3--8--2) and (3--8--3) for the case $s(\sigma )>0$ and
$s(\tau\sigma )>0$ and Theorems 2.1 and 3.5.1---3.5.3 for the case
$s(\sigma )=0$ or $s(\tau\sigma )=0$, this permits to construct
many topological types of
real Enriques surfaces and describe all theoretical possibilities.

The topological classification of
real Enriques surfaces $Y$ with connected non-orientable $Y(\r )$ was
obtained. For these surfaces, $Y(\r )$ is one of the
following non-orientable connected surfaces
(all these possibilities take place):
$U_k,\ k=1,3,5,7,9$. Here $U_k$ denotes a connected non-orientable
surface of
the Euler characteristic $1-k$;
its $2$-sheeted unramified orientable covering is a connected
orientable surface $T_k$ of the
genus $k$. For all these surfaces the invariants
$r(\theta )=a(\theta )$. Thus, by
Theorem 3.6.1, for these surfaces $_2Br(Y)\cong \z/2$. Besides,
we can remove the condition
$r(\theta )=a(\theta )$ from the formulation of Theorem 3.6.1.

Also, real Enriques surfaces $Y$ are constructed with
$Y(\r )=U_k\amalg U_0$, where $k=2,4,6,8,10$ (at first, the type
$Y(\r )=U_{10}\amalg U_0$ was constructed by R. Silhol \cite{Si})
and one of invariants $s(\sigma ), s(\tau\sigma )$ equal to zero. By
this construction, using Theorems 2.1 and 3.5.1--3.5.2, for these
surfaces we have $b(Y)=2$.
Thus, by Theorem 3.8.5, for these Enriques surfaces,
Hochschild--Serre spectral sequence
degenerates and $\dim~_2Br(Y)=3$. But it is not known if
the homomorphism (0--1) epimorphic for these cases.

Enriques surfaces with 6 real connected components were
constructed. For example, there exists a real Enriques surface $Y$ with
$Y(\r )= U_1\amalg 5T_0$. By Theorem 3.7.1, it is the maximum number of
real connected components for real Enriques surfaces.
For this surface $s(\sigma )>0$ and $s(\tau\sigma )>0$.
Thus, by Theorem 3.8.3, $\dim~_2Br(Y)=11$.

Enriques surfaces with 4 real connected non-orientable components were
constructed. For example, there exists a real Enriques surface $Y$ with
$Y(\r )= 2U_2\amalg 2U_0$. By Theorem 3.7.1, it is the maximum number of
real connected non-orientable components for real Enriques surfaces.

We send the reader to the
paper \cite{N6} for further examples.

At last, we want to mention the most interesting problem (from our
point of view) connected with real Enriques surfaces: To construct an
example of a real Enriques surface such that the homomorphism
(0--1) is not an epimorphism. By Theorem 3.8.3 (from \cite{N5}), one can
construct this example only if one of
$s(\sigma ),\ s(\tau\sigma )$ is zero.

\enddocument


\newpage

\Refs

\widestnumber\key{CT-P}

\ref
\key A
\book Algebraic surfaces \ed I. R. Shafarevich
\publ Proc. Steklov Math. Inst. Vol 75 \yr 1965
\transl\nofrills English transl. by A.M.S.
\yr 1969
\endref


\ref
\key Br \by G. E. Bredon \book Introduction to compact
transformation groups
\publ Academic Press, New York and London \yr 1972
\endref


\ref
\key CT-P \by J.-L. Colliot-Th\'el\`ene and R. Parimala
 \paper Real components of algebraic varieties and \'etale cohomology
\jour Invent. math.
\vol  101 \yr 1990 \pages 81--99
\endref

\ref
\key C-D \by F. R. Cossec and I. Dolgachev \paper Enriques
surfaces.I
\inbook Progress in Mathematics \vol 76. \publ Birkh\"auser \yr 1989
\endref


\ref
\key Ha1 \by V. M. Harlamov \paper
Topological types of nonsingular surfaces of degree 4 in $RP^3$
\jour Funktcional. Anal. i Prilozhen.
\vol  10 \yr 1976 \pages 55--68
\transl\nofrills English transl. in \jour Functional Anal. Appl.
\endref

\ref
\key Ha2 \bysame \paper
Real algebraic surfaces
\inbook Proceedings of the International
Congress of Mathematicians, Helsinki, I
\yr 1978 \pages 421-428
\endref

\ref
\key Ho \by E. Horikawa \paper
On the periods of Enriques surfaces, I, II
\jour Math. Ann.
\vol  234 and 235\yr 1978 \pages 73--108 and 217--246
\endref


\ref
\key Kr \by V. A. Krasnov \paper Harnack-Thom inequalities
for mappings of real algebraic varieties \jour Izv. Akad. Nauk SSSR
Ser. Mat. \vol 47 \yr 1983 \pages 268--297
\transl\nofrills English transl. in  \jour Math. USSR Izv.
\vol 22 \yr 1984  \pages 247--275
\endref

\ref
\key Ku \by Vik. S. Kulikov
\paper Degenerations of K3-surfaces and Enriques surfaces
\jour Izv. Akad. Nauk SSSR Ser. Mat.
\vol  41 \yr 1977 \pages 1008--1042
\transl\nofrills English transl. in  \jour Math. USSR Izv.
\vol 11 \yr 1978  \pages 957--989
\endref

\ref
\key Ma \by Y. I. Manin
\paper Le groupe de Brauer-Grothendieck en G\'eometrie
diophantienne
\inbook Actes du Congr\`es Intern. Math. Nice (1970) \vol 1
\yr 1971 \pages 401-411 \publ Gauthier-Villars
\publaddr Paris
\endref

\ref
\key Mi \by J. Milne \book \'Etale cohomology
\yr 1980 \publ Princeton Univ. Press
\endref

\ref
\key M-N \by Sh. Mukai and Yu. Namikawa
\paper Automorphisms of Enriques surfaces which act trivially on
the cohomology groups
\jour Invent. math.
\vol  77 \yr 1984 \pages 383--397
\endref


\ref
\key N0 \by V. V. Nikulin
\paper
On Kummer surfaces
\jour Izv. Akad. Nauk SSSR Ser. Mat.
\vol  39 \yr 1975 \pages 278--293
\transl\nofrills English transl. in
\jour Math. USSR Izv.
\vol 9 \yr 1976  \pages 261--275
\endref



\ref
\key N1 \bysame \paper
Finite groups of automorphisms of K\"ahlerian
surfaces of type K3
\jour Trudi Moscow Mat. Ob-va
\vol 38 \yr 1979 \pages 75--138
\transl\nofrills English transl. in
\jour Moscow Math. Soc. \vol 38 \yr 1980
\endref

\ref
\key N2 \bysame \paper
On the quotient groups of the automorphism groups of
hyperbolic forms by the subgroups generated by 2-reflections,
Algebraic-geometric applications
 \jour Current Problems in Math. Akad. Nauk SSSR, Vsesoyuz. Inst.
Nauchn. i Tekhn. Informatsii, Moscow
\vol  18 \yr 1981 \pages  3-114
\transl\nofrills English transl. in
\jour J. Soviet Math.
\vol 22 \yr 1983 \pages 1401-1476
\endref

\ref
\key N3 \bysame \paper
Integral symmetric bilinear forms and some of their geometric
applications
\jour Izv. Akad. Nauk SSSR Ser. Mat.
\vol  43 \yr 1979 \pages 111--177
\transl\nofrills English transl. in
\jour Math. USSR Izv.
\vol 14 \yr 1980  \pages 103--167
\endref

 \ref
\key N4 \bysame \paper
Involutions of integral quadratic forms and their application to
real algebraic geometry
\jour Izv. Akad. Nauk SSSR Ser. Mat.
\vol  47 \yr 1983 \pages 109--188
\transl\nofrills English transl.in
 \jour Math. USSR Izv.
\vol 22 \yr 1984
\endref

\ref
\key N5 \bysame \paper
Lectures on the Brauer group of real algebraic surfaces
\jour Preprint of University of Notre Dame,
College of science, Dept. of Math. \vol \# 179  \yr January, 1992
\endref

\ref
\key N6 \bysame \paper
On the topological classification of real Enriques surfaces.
\jour Preprint \yr 1993
\endref

\ref
\key N-S \by V. V. Nikulin and R. Sujatha
\paper On Brauer groups of real Enriques surfaces
\jour Preprint IHES /M/91/33 \yr 1991
\endref

\ref \key P\u S-\u S
\by I. I. Pjatecki\u i-\u Sapiro and I. R. \u Safarevi\u c
\paper
A Torelli theorem for algebraic surfaces of type K3
\jour Izv. Akad. Nauk SSSR Ser. Mat.
\vol  35 \yr 1971 \pages 530--572
\transl\nofrills English transl. in
\jour Math. USSR Izv.
\vol 5 \yr 1971
\endref

\ref
\key Si
\by R. Silhol
\paper Real algebraic surfaces
\inbook Lecture Notes Math. \vol 1392 \yr 1989 \pages 215
\endref

\ref
\key Sp
\by E. H. Spanier \book Algebraic topology
\publ McGraw-Hill Book Company \yr 1966
\endref

\ref
\key Su \by R. Sujatha \paper
Witt groups of real projective surfaces
\jour Math. Ann.
\vol  28 \yr 1990 \pages 89--101
\endref




\endRefs
\bye